\documentclass[twocolumn,final]{svjour3}

\usepackage{csquotes}

\usepackage{amsmath}
\usepackage{amssymb}
\usepackage{mathtools}
\usepackage{tensor}
\usepackage{colonequals}

\usepackage{graphicx}
\usepackage{array} 
\usepackage{booktabs} 
\usepackage{multirow}
\usepackage{caption}
\usepackage{subcaption}

\usepackage{algorithmic}

\usepackage{tikz}
\usetikzlibrary{arrows,shapes,plotmarks}
\tikzset{>=stealth'} 
\tikzstyle{graphnode} = 
   [circle,draw=black,minimum size=22pt,text centered,text
     width=22pt,inner sep=0pt] 
\tikzstyle{var}   =[graphnode,fill=white]
\tikzstyle{obs}   =[graphnode,fill=black,text=white]
\tikzstyle{fac}   =[rectangle,draw=black,fill=black!25,minimum size=5pt]
\tikzstyle{facprior} =[rectangle,draw=black,fill=black,text=white,minimum size=5pt]
\tikzstyle{edge}  =[draw=white,double=black,thick,-]
\tikzstyle{prior} =[rectangle, draw=black, fill=black, minimum size=
5pt, inner sep=0pt]
\tikzstyle{dirprior} = [circle, draw=black, fill=black, minimum
size=5pt, inner sep=0pt]

\usepackage{pgfplots}
\pgfplotsset{compat=newest}
\pgfplotsset{plot coordinates/math parser=false}     

\newlength\figureheight
\newlength\figurewidth

\usepackage[ruled,vlined]{algorithm2e}

\definecolor{lred}{RGB}{200,0,0}
\definecolor{dred}{RGB}{130,0,0} \definecolor{dblu}{RGB}{0,0,130}
\definecolor{dgre}{RGB}{0,130,0} \definecolor{dgra}{RGB}{50,50,50}
\definecolor{mgra}{RGB}{100,100,100}
\definecolor{lgra}{RGB}{220,220,220}
\definecolor{MPG}{RGB}{000,125,122}
\definecolor{ora}{HTML}{FF9933}

\definecolor{AMPurple}{HTML}{663366}
\definecolor{Burgundy}{HTML}{993333}
\definecolor{Coffee}{HTML}{7B6049}
\definecolor{ForestGreen}{HTML}{005826}
\definecolor{Lavender}{HTML}{6E6AB1}
\definecolor{PSLightBlue}{HTML}{7DA7D9}

\usepackage{tgheros}

\usepackage{slantsc}
\usepackage{microtype}

\usepackage[setpagesize=false,
            colorlinks=true,
            linkcolor=black,
            citecolor=black,
            urlcolor=black,
            hyperfootnotes=false]{hyperref}
\usepackage{url}

\usepackage{pdflscape}

\usepackage{microtype}

\usepackage{nicefrac}

\newcommand{\matlab}{\texttt{MATLAB}}

\renewcommand{\Re}{\mathbb{R}} 

\newcommand{\Exp}{\mathbb{E}} 
\newcommand{\cov}{\operatorname{cov}} 

\newcommand{\N}{\mathcal{N}} 

\newcommand{\Trans}{^{\intercal}} 
\newcommand{\mat}{\vec}
\newcommand{\Id}{\vec{I}} 

\newcommand{\g}{\,|\,} 
\newcommand{\ce}{\colonequals} 
\newcommand{\ec}{\equalscolon} 
\newcommand*{\diff}{\mathop{}\!\mathrm{d}}

\newcommand{\abs}[1]{\lvert #1\rvert} 
\newcommand{\norm}[1]{\lvert\lvert #1\rvert\rvert}
\newcommand{\ind}[1]{\mathbb{I}_{#1}} 

\newcommand{\q}{\quad}

\newcommand{\ssq}{\sigma^2}

\newcommand{\tn}{{t_n}}
\newcommand{\tnmo}{{t_{n-1}}}
\newcommand{\tnpo}{{t_{n+1}}}

\newcommand{\tmo}{{t_{-1}}}

\newcommand{\mesh}{\mathrm{\Delta}}

\newcommand{\matexp}[2]{e^{\mat{F}\Delta^{#1}_{#2}}}

\pgfkeys{/pgfplots/mystyle/.style={
  semithick,
  tick style={major tick length=4pt,semithick,gray},
  xtick align = inside,
  ytick align = inside,
  xlabel near ticks,
  ylabel near ticks
  }}

\newif\iffinal 

\usepackage{color}

\iffinal
  \newcommand{\simo}[1]{}
  \newcommand{\michael}[1]{}
  \newcommand{\hennig}[1]{}
  \newcommand{\yell}[1]{}
  \newcommand{\editor}[1]{}
  \newcommand{\revone}[1]{}
  \newcommand{\revtwo}[1]{}
\else
  \newcommand{\simo}[1]{{\color{dblu} \textbf{[Simo wonders: #1]}}}
  \newcommand{\michael}[1]{{\color{dred} \textbf{[Michi muses: #1]}}}
  \newcommand{\hennig}[1]{{\color{MPG} \textbf{[Philipp says: #1]}}}
  \newcommand{\yell}[1]{{\color{ora} \textbf{[#1]}}}
  \newcommand{\editor}[1]{{\color{ora} \textbf{[Editor: #1]}}}
  \newcommand{\revone}[1]{{\color{ora} \textbf{[Reviewer 01: #1]}}}
  \newcommand{\revtwo}[1]{{\color{ora} \textbf{[Reviewer 02: #1]}}}
\fi

\newcommand{\editone}[1]{\textcolor{blue}{#1}}
\newcommand{\edit}[1]{\textcolor{red}{#1}}
\usepackage[normalem]{ulem} 
\newcommand{\del}[1]{\textcolor{blue}{\sout{#1}}}

\renewcommand{\editone}[1]{#1}
\renewcommand{\edit}[1]{#1}
\renewcommand{\del}[1]{}

\newcommand{\repo}{\footnote{\url{https://pn.is.tuebingen.mpg.de/code/pfos}}}

\title{A probabilistic model for the numerical solution of initial value problems}

\titlerunning{Stat Comput}
\authorrunning{Stat Comput}

\author{Michael Schober%
  \and
  Simo S\"arkk\"a%
  \and
  Philipp Hennig
}

\institute{Max Planck Institute for Intelligent Systems,\\
  72076 T\"ubingen, Germany
           \and
  Aalto University,\\
  Finland
}

\date{Received: date / Accepted: date}

\begin{document}
\maketitle

\begin{abstract}
We study connections between ordinary differential equation (ODE)
solvers and probabilistic regression methods in statistics.
We provide a new view of probabilistic ODE solvers as active
inference agents operating on stochastic differential equation
models that estimate the unknown initial value problem (IVP)
solution from approximate observations of the \edit{solution derivative},
as provided by the ODE dynamics.
Adding to this picture, we show that several multistep methods
of Nordsieck form can be recasted as Kalman filtering on
$q$-times integrated Wiener processes.
Doing so provides a family of IVP solvers that return a Gaussian
posterior measure, rather than a point estimate.
We show that some such methods have low computational overhead,
nontrivial convergence order, and that the
posterior has a calibrated concentration rate.
Additionally, we suggest a step size adaptation algorithm
which completes the proposed method to a practically useful
\edit{implementation}, which we experimentally evaluate using a
representative set of standard codes \edit{in} the DETEST benchmark set.
\end{abstract}

\keywords{
Initial value problems \and
Nordsieck methods \and
Runge--Kutta methods \and
filtering \and
Gaussian processes \and
Markov processes \and
probabilistic numerics%
}

\CRclass{60H30 \and 62M05 \and 65C20 \and 65L05 \and 65L06}



\section{Introduction}
\label{sec:introduction}

Numerical algorithms estimate intractable quantities from tractable ones.
It has been pointed out repeatedly
\cite{poincare1896,diaconis88:_bayes,ohagan92:_some_bayes_numer_analy}
that this
process is structurally similar to statistical inference, where the
tractable computations play the role of data in statistics, and the
intractable quantities relate to latent, inferred quantities.
In recent years, the search for numerical algorithms which
return probability distributions over the solution for a given numerical problem
has become an active area of research \cite{hennig15probabilistic}.
Several models and methods have been proposed
for the solution of initial value problems (IVPs) 
\cite{skilling1991bayesian,o.13:_bayes_uncer_quant_differ_equat,%
schober2014nips,conrad_probability_2015,KerstingHennigUAI2016,teymur2016probabilistic}.
\editone{However, these \edit{probabilistic} algorithms have no immediate
connection to the \edit{extensive} literature on this task
in numerical analysis. Most importantly, such
inference algorithms do not come with
convergence analysis out of the box.}
The methods in
\cite{o.13:_bayes_uncer_quant_differ_equat,conrad_probability_2015,teymur2016probabilistic}
have convergence results, but
\editone{their respective implementations are based on
sampling schemes and, thus, do not offer guarantees for
individual runs.}
The methods in \cite{schober2014nips,KerstingHennigUAI2016}
\editone{offer a deterministic execution and an
analytical guarantee for the first step, but
we will show that this guarantee is lacking
for the whole integration domain.}

In this \edit{paper}, we present a class of probabilistic solvers which
combine properties of the standard and the probabilistic \edit{algorithms}.
\editone{We formulate desiderata that users might have for
a probabilistic numerical algorithm. We present one construction
that fulfills these desiderata and we provide a \matlab{} code\repo{}
\edit{which} we compare empirically against other available codes.
The construction uses the algebra of Gaussian inference to provide
a Gaussian posterior distribution over the solution of an IVP.}
\edit{In particular, we show that the posterior mean
can be understood as a multistep method in Nordsieck representation
and, thus, analytical results about these methods carry over to
the present algorithm.}
Additionally, we propose to interpret the posterior covariance as a
\emph{measure of uncertainty} or \emph{error estimator}
and \editone{argue} that this interpretation \editone{can be} analytically justified.
In the context of a larger pipeline of empirical studies and
numerical computations, the framework of probability modeling provides
a common language to analyze the epistemic confidence in its result
\cite{cockayne2017bayesian}.
\editone{In the framework of Cockayne et al.~\cite{cockayne2017bayesian}, the code
provides approximate Bayesian uncertainty quantification
\cite{sullivan2015introduction}
at low computational overhead and almost complete backwards
compatibility to the \matlab{} IVP solver suite.}


\subsection{Problem description}

We study the problem of finding a real-valued curve
$y:\mathbb{T}\to\Re$ over \edit{an interval $\mathbb{T} = [t_0, T]$} such that
\begin{subequations}
\begin{align}
y'(t) &= \frac{\diff y}{\diff t}(t) = f(t,y(t)) \text{~(``the ODE''),}\\
\intertext{and}
y(t_0) &= y_0 \text{~(``the initial value'')}, \label{eq:ivp}
\end{align}
\label{eq:ivp-problem}  
\end{subequations}
with $f$ Lipschitz continuous \edit{with constant $L$ in the second argument}
and sufficiently many times differentiable
in its second argument.
\editone{
\edit{U}sers might be
interested in approximations to $y$ on either a predefined mesh
$\mesh_S \subset \mathbb{T}$ or an automatically selected mesh
$\mesh \subset \mathbb{T}$ of finitely many intermediate function
values.
}
The derivations will
be presented with a scalar-valued problem, but the results carry over
to the \edit{multivariate} case.

\editone{
  \edit{IVPs} are a particularly deeply studied class
  of ODE-related tasks. Part of their significance is due
  to the Picard--Lindel\"of theorem which guarantees
  local unique existence of solutions.
  As a consequence, IVPs lend
  themselves to be solved by so-called \emph{step-by-step methods},
  where the solution is advanced iteratively on expanding meshes
  $\mesh_{n+1} \ce \edit{(\{t_0, \dotsc, \tn\} \cup \{\tnpo\})} \supset \mesh_n$.
}%
\edit{
  A typical mesh is generated by choosing a \emph{step size} $h \in (0, T - t_0]$
  and by setting $t_n \ce t_0 + hn$.
  Alternatively, $h_n$ may vary per step and $t_n = t_0 + \sum_{i=1}^n h_i$.
}

\editone{
  To construct a probabilistic numerical method, we define the
  following list of desiderata that an algorithm should fulfill.
  These properties will be defined and motivated in turn below.
  \begin{description}
    \item[\textbf{Probabilistic inference.}]
      \edit{The} computations should
      be operations on probability distributions.
    \item[\textbf{Global definition.}]
      The probabilistic model should not depend on the
      discretization mesh.
    \edit{
    \item[\textbf{Deterministic execution.}]
      \edit{When} run several times on the same problem, the algorithm produces the same output each time.
    \item[\textbf{Analytic guarantees.}]
      The algorithm's output should have
      desirable analytic properties.
    }
    \item[\textbf{Online execution.}]
      The algorithm execution can be
      extended indefinitely when required.
    \item[\textbf{Speed.}]
      The execution time should not be prohibitively slow.
    \item[\textbf{Problem \edit{adaptiveness}.}]
      \edit{The algorithm should automatically
        adapt parameters to problem and accuracy requirements.}
  \end{description}
}

Throughout this \edit{paper}, we will use zero-based indexing for vectors
and matrices such that
a $d$-di\-men\-sion\-al vector $\vec{v}$ is written as
$\vec{v} = (v_0, \dotsc, v_{d-1})\Trans$ and the $d$ canonical basis
vectors are $\vec{e}_0, \dotsc, \vec{e}_{d-1}$.






\section{From classical to probabilistic numerical algorithms}
\label{sec:background-statistics}

\editone{
In this section, we explain and motivate the first two
items from our list of desiderata in turn---\emph{probabilistic inference}
and \emph{global definition}.
}

\editone{
  On a high-level view, numerical algorithms can be
  described as combinations of
  \emph{tractable approximating function classes}
  and \emph{computation strategies} for \emph{informative values}.
  Analyses of numerical methods show to what level
  the approximations can converge to the true problem
  solution and how fast the computation strategies can
  be carried out.
  This is structurally very similar to problems in
  statistics where unknown quantities need to be
  related to \emph{approximating function classes}
  via \emph{observable informative values}.
  In particular, finding a function
  $Y = (Y_t)_{t \in \mathbb{T}}$ given a
  collection of information $z_n,\;n = 0, \dotsc N$
  about $Y_\tn$ at times $\tn$ is studied in
  \emph{regression analysis} in statistics.
  \edit{
  In that context, the unknown function is often treated
  as a stochastic process and the approximating function
  is obtained by conditioning it on the measurements.
  }%
  Consequently, this paper treats the problem
  of finding an approximate solution $Y = (Y_t)_{t \in \mathbb{T}}$
  to the true unknown solution $y(t)$ as a \edit{statistical} regression problem
  \edit{on a stochastic process}.
}


\del{
Many considerations must be taken into account when choosing
an appropriate model. In terms of statistical power, the
optimal model is the true generative process,
if it is known.
In this work, the model is selected based on desiderata
derived from applications and competing algorithms.
}
Accepting the probabilistic approach as a framework
for plausible reasoning
\cite{jeffreys1969theory,cox1946probability,hennig15probabilistic},
we require a \edit{\emph{probability measure} or \emph{law} $P_Y$} over
the numerical solution \edit{$Y_t$}.
\editone{
  The computations necessary for the construction of $P_Y$
  should be interpretable as (approximate) probabilistic
  inference.
  When such an interpretation is admissible, we call the
  resulting algorithm a \emph{probabilistic numerical
  method (PNM)} for the purposes of this paper.
  A more rigorous definition has been given by
  Cockayne et al.~\cite{cockayne2017bayesian}.
  The motivation behind this requirement is that there
  should not be an analysis gap between statistical
  and numerical computations.
}
This is particularly
beneficial, when the differential equation solver
is embedded in a longer chain of computations
\editone{
\cite{cockayne2017bayesian}.
In principle, this should allow to build fine-tuned
methods adopting to sources of data uncertainty and
computational approximation during runtime and provide
richer feedback of approximation quality as recently empirically validated
by \cite{schober2014probabilistic,hauberg2015random}.
}

Let $z_{[n]} \ce \{ z_k\g k \le n\}$ be the set of
collected data up to and including step $n$.
Given a \emph{prior \editone{law}} \editone{$P_Y$} over the
space of solutions and a \emph{likelihood function}
$P(z_n\g Y_\tn)$ relating the
\edit{value of the process $Y_\tn$}
to collected data, Bayes' theorem leads to the (predictive) posterior \edit{measure}
\begin{equation}
\editone{P_{Y\g z_{[n]}}} =
  \frac{P(z_{[n]}\g Y)\,\editone{P_Y}}
       {\int P(z_{[n]}\g Y)\,\editone{\diff} \editone{P_Y}}.
\label{eq:bayes}
\end{equation}
\edit{
where $P(z_{[n]}\g Y) = \prod_{k \leq n} P(z_k\g Y_{t_k})$. 
Rigorously, the above expression is valid only for finite
collections of values of $Y_t$, in which case the corresponding
probability measures $P_Y$ are typically represented by their densities,
but as the finite-dimensional distributions define the full measure,
we use this slight abuse of notation here.
We denote the posterior distributions (typically densities)
of point values of $Y$ as $P(Y_t \g z_{[n]})$.
}


\edit{
The objects of interest in this case will be a
stochastic process, so some measure-theoretic
restrictions apply to Equation~\eqref{eq:bayes}
\cite{stuart_2010}, \cite[\textsection 7.3]{gine2015mathematical}.
}

We\del{rather} propose to think about the probabilistic
framework\del{only} as a more informative output information
than the point estimates returned by classical
numerical algorithms (see also~\cite{hennig15probabilistic}).

\editone{
  \edit{Furthermore, a} probabilistic IVP solver shall be called
  \emph{globally defined} on its input domain $\mathbb{T}$,
  if its probabilistic interpretation does not depend
  on the discretization mesh $\mesh$.
  PNMs satisfying this property provide two benefits.
  Users may evaluate the (predictive posterior) \edit{distribution $P(Y_t\g z_{[n]})$}
  for any $t \in \mathbb{T}$. In particular, users
  may evaluate \edit{$P(Y_t\g z_{[n]})$} for $t \notin \edit{\mesh}$.
  Thus, users may request \edit{$P(Y_{t_s}\g z_{[n]}), t_s \in \mesh_S$} and
  the support of a user-defined mesh $\mesh_S$ is not a
  separate requirement. Secondly, this implies that the
  inference can be paused and continued after every
  expansion from $\mesh_n \mapsto \mesh_{n+1}$.
  In principle, this also enables iterative refinement
  of the solution quality based on its prediction
  uncertainty.
}

\editone{
  Table~\ref{tab:pnm-ode-solvers} lists PNM ODE solvers that
  have been proposed in the literature. A $\checkmark$
  indicates that the method satisfies a given property,
  a $\times$ indicates that a method does not satisfy a
  given property, and a $\thickapprox$ indicates that a
  property holds with some restrictions.
  The listing shows that almost all methods proposed so
  far are globally defined. Furthermore, we see that
  the definition is independent of a method being
  sampling-based or not.
  The method proposed by Conrad et al.~\cite{conrad_probability_2015}
  is a generative process on sub-intervals
  $[\tn, \tnpo] \subset \mathbb{T}$ based on a
  numerical discretization. It is easy to construct
  two different meshes $\mesh_n, \mesh_n'$ that define
  different distributions for \edit{$Y_t$}
  in the case   of $y' = \lambda y$
  and a general argument can be made from this example.
  In Teymur et al.~\cite{teymur2016probabilistic}, the predictive
  posterior is only defined on the discretization mesh.
  This defect is not for lack of definition, but 
  a consequence of the underlying numerical
  method the probabilistic algorithm is built upon.
  Since the method is defined on a windowed data frame,
  it is easy to construct a mesh such that the prediction
  \edit{$Y_t$} at time $t$ will be different depending on the
  window $[t_{n-i}, \dots, t_{n+j}] \ni t$ is
  considered to be part of.
}

\editone{
  The analysis in Schober et al.~\cite{schober2014nips}
  proposes two main modes of operation: naive chaining and probabilistic
  continuation. Naive chaining is not a globally defined
  method since mesh points $t_n$ are part of adjacent
  Runge--Kutta blocks and the corresponding predictive
  posterior \edit{distribution $P(Y_\tn\g z_{[n]})$} is different for these two
  blocks. Probabilistic continuation is globally
  defined, but there has been no convergence theory
  for this mode yet. This paper fills this gap.
}

\begin{table}[t]
\caption{Properties of existing PNM ODE solvers}
\label{tab:pnm-ode-solvers}
\begin{center}
\editone{
\begin{tabular}{lccc}
  \toprule
  Method & glob.\ def.? & determ.? & guarantees?\\
  \midrule
  Skilling \cite{skilling1991bayesian} & $\checkmark$ & $\times$ & $\times$ \\
  Chkrebtii et al.~\cite{o.13:_bayes_uncer_quant_differ_equat} & $\checkmark$ & $\times$ & $\thickapprox$ \\
  Schober et al.~\cite{schober2014nips} & $\thickapprox$ & $\checkmark$ & $\thickapprox$ \\
  Conrad et al.~\cite{conrad_probability_2015} & $\times$ & $\times$ & $\thickapprox$ \\
  Kersting & $\checkmark$ & $\checkmark$ & $\times$ \\
  \& Hennig \cite{KerstingHennigUAI2016}\\
  Teymur et al.~\cite{teymur2016probabilistic} & $\times$ & $\times$ & $\thickapprox$ \\
  \textbf{PFOS (this paper)} & $\checkmark$ & $\checkmark$ & $\checkmark$ \\
  \bottomrule
\end{tabular}
}
\end{center}
\end{table}


\subsection{State-space models for Gauss--Markov processes}
\label{sec:state-space}

\edit{
Our approximate model of the true solution $y(t)$ is
a vector $\vec{x}(t) = (y^{(0)}(t), \dotsc, y^{(q)}(t))\Trans$
where $y^{(i)}(t)$ is the true $i$-th derivative of $y(t)$
at time $t$.
We represent the prior uncertainty about $\vec{x}(t)$ by
the distribution $P(\vec{X}_t)$ of the random variable $\vec{X}_t$---or
more generally as the measure or the law $P_{\vec{X}}$ of the stochastic process
$\vec{X}$---which is then conditioned on the observed values.
}

The \edit{prior} model\edit{, which has also been considered}
in Schober et al.~\cite{schober2014nips}, belongs to the class of
Gauss--Markov \editone{processes. Models of this class can often be written} as
a \emph{linear time-invariant (LTI)}
\emph{stochastic differential equation (SDE)} of the form
\begin{equation}
\diff\vec{X}_{\edit{t}} = \mat{F} \, \vec{X}_{\edit{t}} \diff t + \vec{L} \, \diff W_{\edit{t}},
\label{eq:sde}
\end{equation}
where \editone{$\vec{X}_t$} is the so-called \emph{state} of the model,
$\mat{F} \in \Re^{(q+1) \times (q+1)}$ is the state feedback matrix
and $\vec{L} \in \Re^{(q+1)}$ is the diffusion matrix
of the system. $\diff \edit{W_t}$ is the increment of a Wiener process with intensity
$\sigma^2$, that is, $\diff \edit{W_t} \sim \N(0, \sigma^2 \diff t)$.

\edit{
Here, we consider models where $\vec{L}$ is the last
standard basis vector $\vec{e}_q$ 
and $\mat{F} = \mat{U}_{q+1} + \vec{e}_q\vec{f}\Trans$
is a (transposed) companion matrix. Here, $\mat{U}_{q+1}$ denotes
the upper shift matrix and the row vector $\vec{f}\Trans$ 
contains the coefficients in the last row of $\mat{F}$.
In this case, the vector-valued process
$\vec{X}_t = (X_{t,0}, \dotsc, X_{t,q})\Trans$ obtains the
interpretation $\vec{X}_t = (Y_t, Y_t', \dotsc, Y^{(q)}_t)\Trans$,
because the form of $\mat{F}$ and $\vec{L}$ implies that
the realizations of $Y_t$ are $q$-times continuously
differentiable on $\Re$.
Later, we will also consider scaled systems $\vec{\tilde{X}}_t = \mat{B} \vec{X}_t$
with an invertible linear transformation $\mat{B}$.
In this case, we denote by $\mat{H}_i$ the matrix
that projects onto the $i$-th derivative
$Y^{(i)}_t = \mat{H}_i \vec{\tilde{X}}_t \ce e_i\Trans \mat{B}^{-1} \vec{\tilde{X}}_t$.
}%
Two particular models of this type
are the $q$-times \emph{integrated Wiener process} (IWP($q$)) and the
continuous auto-regressive processes of order $q$.
Detailed introductions can be found, for example, in
\cite{karatzas1991brownian,oksendal2003stochastic,sarkka2006thesis}.
SDEs can also be seen as path-space representations
of more general temporal Gaussian processes
arising in machine learning models \cite{sarkka2013spm}.

\editone{
  Models of the form \eqref{eq:sde}
  are also related to
  non-parametric spline regression models
  \cite{wahba1990spline} which often have a natural
  interpretation in frequentist analysis \cite{kimeldorf1970correspondence}.
  Conceptually, these models are a compromise between
  globally defined parametric models, which might be
  too restrictive to achieve convergence, and local
  parametric models, which might be too expressive to
  be captured by a globally defined probability distribution.
  Models of this type have been studied in the
  literature \cite{Loscalzo1967Spline,andria1973integration},
  but the presentation here starts from
  other principles.
}

Conditioning on (random) initial conditions $\vec{X}_{t_*}$ at a starting time
$t_*$ of the process, the solution of Equation~\eqref{eq:sde} has the analytic form
\begin{equation}
  \vec{X}_t = \matexp{t}{t_*}\vec{X}_{t_*}
  + \int_{t_*}^t \matexp{t}{\tau} \, \vec{L} \, \diff W(\tau),
  \label{eq:sde-sol-general}
\end{equation}
\editone{
where
$\matexp{t}{t'} \ce
 \sum_{k=0}^\infty [\mat{F}\edit{\Delta^t_{t'}}]^k [k!]^{-1}$
is the matrix exponential of $\mat{F}\Delta^t_{t'}$
\edit{and $\Delta^t_s \ce t - s$.}
}

If $\vec{X}_{t_*} \sim\N(\vec{m}_*, \mat{C}_*)$, then the
distribution of $\vec{X}_t$ remains Gaussian for all $t$ by linearity
and its statistics can be computed explicitly
\cite{grewal2001kalman,sarkka2006thesis} via
\begin{equation}
\begin{aligned}
\vec{m}_t \ce &\Exp(\vec{X}_t) = \matexp{t}{t_*} \vec{m}_*\\
\cov(\vec{X}_t,\vec{X}_{t'}) = 
  &\matexp{t}{t_*}\mat{C}_*(\matexp{t'}{t_*})\Trans\\
   &+ \underbrace{\int_{t_*}^{\min(t,t')}\matexp{t}{\tau}\vec{L}\ssq\vec{L}\Trans(\matexp{t'}{\tau})\Trans d\tau}%
   _{\ec \mat{Q}_{t_*}(t,t')}.
\end{aligned}
\label{eq:sde-sol-mu-cov}
\end{equation}
For practical purposes, only the covariance matrix $\mat{C}_t = \cov(\vec{X}_t,\vec{X}_t)$
of the states at a single time $t$ is needed.

\edit{
The choice of prior measure $P_{\vec{X}}$ in Equation~\eqref{eq:sde} can be 
interpreted as a \emph{prior assumption} or \emph{belief}
encoded in the algorithm, in the sense that the algorithm amounts
to an autonomous agent. We emphasize that, if one adopts this view,
then the results reported in later sections amount to an external
analysis of the effects of these assumptions. That is, we will show
that if the agent ``makes'' these assumptions, they give rise to a
posterior distribution with certain desirable properties.
By contrast, one could also take a more restrictive standpoint
internal to the algorithm, and state that the proposed method
works well if the true solution $\vec{x}(t)$ is indeed a sample from $P_{\vec{X}}$.
This is expressively \emph{not} our viewpoint here;
and it would be a flawed argument, too, given that in practice,
$\vec{x}(t)$ is defined through the ODE, thus evidently not a sample from any stochastic process.
}

Denote by $\mat{A}(h) \ce \matexp{t + h}{t}$ the \emph{discrete transition matrix} of
step size $h$ and $\mat{Q}(h) \ce \mat{Q}_t(t+h,t+h)$ the \emph{discrete diffusion matrix}
of step size $h$, respectively. For LTI SDE systems, $\mat{A}(h)$ and $\mat{Q}(h)$
fulfill matrix-valued differential equations which can be solved analytically via
matrix fraction decomposition \cite{grewal2001kalman,sarkka2006thesis}. \edit{If we define}
\begin{equation}
\mat{\Phi}(h) = \begin{pmatrix} \mat{\Phi}_{11}(h) & \mat{\Phi}_{12}(h) \\
                                \mat{\Phi}_{21}(h) & \mat{\Phi}_{22}(h) \end{pmatrix}
  \ce \exp\left\{\begin{pmatrix} \mat{F} & \ssq\vec{L}\vec{L}\Trans \\ 
                                 \mat{0} & -\mat{F}\Trans \end{pmatrix} h \right\},
\end{equation}
\edit{then, the matrices} $\mat{A}(h)$ and $\mat{Q}(h)$ are given by
\begin{equation}
\mat{A}(h) = \exp( \mat{F} h ),\quad \mat{Q}(h) = \mat{\Phi}_{12}(h) \mat{\Phi}_{22}^{-1}(h).
\label{eq:trans-and-diffusion}
\end{equation}
\edit{Above, }$\mat{\Phi}_{22}^{-1}(h)$ can be computed ef\-fi\-cient\-ly:
from the two properties of the matrix exponential, $\exp(\mat{X})^{-1} = \exp(-\mat{X})$
and $\exp(\mat{X}\Trans) = \exp(\mat{X})\Trans$, it follows that
$\mat{\Phi}_{22}^{-1}(h) = \mat{A}(h)\Trans$ and therefore
$\mat{Q}(h) = \mat{\Phi}_{12}(h)\mat{A}(h)\Trans$. In the following,
it will be beneficial to write $\mat{Q}(h)$ as
$\mat{Q}(h)\Trans=\mat{A}(h)\mat{\Phi}_{12}(h)\Trans$, which is valid
since $\mat{Q}(h)$ is symmetric.


\edit{
For the rest of this paper, we will focus on the 
$q$-times integrated Wiener process IWP($q$),
which is defined by
\begin{equation}
\diff \vec{X}_t = \mat{U}_{q+1} \vec{X}_t \diff t + \vec{e}_q \diff W_t.
\label{eq:iwp-q}
\end{equation}
In this case, $\vec{f}\Trans = (0, \dotsc, 0)$ and there is
no feedback from higher states $X_{t,i}$ to lower states $X_{t,j}, i < j$.
In particular, this process is non-stationary and does not
revert to the initial mean $\vec{m}_{t_*}$.
In this system, $\mat{A}(h)$ and $\mat{Q}(h)$ can be be
computed analytically
\begin{equation}
\begin{aligned}
  (\mat{A}(h))_{i,j} &= \ind{i \le j} \frac{h^{\editone{j-i}}}{(j-i)!},\\
  (\mat{Q}(h))_{i,j} &= \ssq \frac{h^{2q+1-i-j}}{(2q+1-i-j)(q-i)!(q-j)!},
\end{aligned}
  \label{eq:iwp-drift-diffusion}
\end{equation}
which can be derived directly from Equation~\eqref{eq:sde-sol-mu-cov}.
}



\subsection{Data generation mechanism}
\label{sec:data-generation}

\editone{
Many problems in statistics assume the existence
of an \edit{externally produced, thus fixed} data set $\{(\tn, z_n)\g\tn \in \edit{\mesh}\}$
and develop appropriate solutions from there.
An analogous concept in numerical algorithms
for solving differential equations would be
to pose a global discretization scheme and to
obtain a solution with other tools from numerical
analysis. Methods of this type are often applied
to boundary value problems (BVPs) and partial
differential equations (PDEs) where the integration
domains need to be specified a priori in any case.
Cockayne et al.~\cite{cockayne2017bayesian} take
this approach by assuming a fixed information
operator $A$.
However, there are cases where the end $T$ of
the integration domain $\mathbb{T}$ cannot be
stated beforehand, when the quantity of interested
depends on a qualitative behavior of the solution.
For example, in modeling of chemical reactions a user
might be interested in the long-term behavior
of the compounds and it is unknown when the reaction
reaches equilibrium.
}

\begin{algorithm*}
\caption{Active probabilistic model}
\label{alg:pseudo}
\begin{algorithmic}[1]
\STATE{Define $t_{-1} \ce t_0,\;z_{-1} \ce y_0$ and probabilistic model \edit{$P_{\vec{X}}$}}
\STATE{Compute $P(\vec{X}_{t_{-1}}\g z_{[-1]}).$}\hfill\COMMENT{Add initial value information}
\FOR{n = 0 \TO \edit{$N$}}
\STATE{Compute $P(\vec{X}_\tn\g z_{[n-1]}) \propto P(\vec{X}_\tn\g\vec{X}_\tnmo) P(\vec{X}_\tnmo \g z_{[n-1]})$}\hfill\COMMENT{\textbf{Predict} $\tn$}
\STATE{\edit{Compute observation model $P(z_n\g\vec{X}_\tn) = \operatorname{observe}(f, P(\vec{X}_\tn\g z_{[n-1]}))$}}\label{alg:line-evaluate}\hfill\COMMENT{Evaluate\editone{/interrogate} model}
\STATE{Compute $P(\vec{X}_\tn\g z_{[n]}) \propto P(z_n\g\vec{X}_\tn) P(\vec{X}_\tn\g z_{[n-1]})$}\hfill\COMMENT{\textbf{Update} information}
\ENDFOR
\RETURN \edit{$\{ P(\vec{X}_\tn\g z_{[n]}),\; n = -1,0,\dotsc,N \}$}
\end{algorithmic}
\end{algorithm*}

\editone{
In contrast, many numerical IVP solvers proceed in
a step-by-step manner.
Having computed a numerical approximation
$P_{Y\g z_{[n]}}$ on the mesh $\mesh_n$,
a \emph{prediction} $y^-_{n+1}$ of $y(\tnpo)$ is
used to \emph{evaluate} $f(\tnpo, y^-_{n+1})$ and
the resulting output $z_{n+1}$ is used to \emph{update}
the approximation $P_{Y\g z_{[n+1]}}$ on the extended mesh $\mesh_{n+1}$.
For example, in a deterministic IVP the data $(t_0, y_0)$ can be
used to construct the observation $z_0 = f(t_0,y_0)$ which satisfies
the probabilistic interpretation of $y'(t_0) \sim \delta(z_0 - y'(t_0))$.
This serves as a corner case for the general situation.
Setting $\tmo \ce t_0$ and $z_{-1} \ce y_0$, it follows
that $y(t_0) \sim \delta(z_{-1} - y(\tmo))$ and the initial
value requires almost no special treatment.
The concept is illustrated in Algorithm~\ref{alg:pseudo} and
can, \edit{in principle,}
be extended indefinitely, \edit{at constant cost per step}. The term 
\emph{Predict-Evaluate-Correct (PEC)} or \emph{Predictor-Corrector}
methods have a more technical meaning in classic textbooks
\cite{hairer87:_solvin_ordin_differ_equat_i,DeuflhardScientific},
but the idea is common to many numerical IVP solvers.
Chkrebtii et al~\cite{o.13:_bayes_uncer_quant_differ_equat}
calls the process of evaluating $f(\tn,y^-_\tn)$ with tentative $y^-_\tn$
to generate $z_n$ a \emph{model interrogation}.
From a statistical perspective, this concept of
active model interrogation is similar
to the sequential analysis of
Wald \cite{wald1973sequential,Owhadi-Scovel-TowardsMachineWald}.
}

\begin{figure}
  \centering

\tikzstyle{smallvar} =
   [circle,draw=black,fill=white,minimum size=22pt,
     text centered]

\begin{tikzpicture}
	\node[text centered, text width=22pt] at (0,2.4) {$\vec{X}_{t_n}$};
	\node[var,minimum width=124pt] at (0,0) (xt) {};

	\node[fac] at (-2.75,0) (ftm1) {} edge (xt);
	\node[var,] (xtm1) at (-3.5,0) {$\vec{X}_\tnmo$} edge (ftm1);

	\node[fac] at (2.75,0) (ftp1) {} edge (xt);
	\node[var] at (3.5,0) (xtp1) {$\vec{X}_\tnpo$} edge (ftp1);
    
    \node[smallvar] (yt0) at (0,1.5) {$\scriptstyle Y^{(0)}_\tn$} edge[dashed] (ftm1);
    \node[smallvar] (yt1) at (0,.25) {$\scriptstyle Y^{(1)}_\tn$} edge[dashed] (ftm1);
    \node[text centered] at (0,-0.5) {$\vdots$};
    \node[smallvar] at (0,-1.5) {$\scriptstyle Y^{(q)}_\tn$} edge[dashed] (ftm1);   
    
    \node[obs] (zn) at (0,-3) {$z_n$};

    \node[fac] (fzn) at (-2.75,-1.5) {} edge (yt1) edge (zn);

    \draw (yt0) edge[->,out=-45,in=45] (zn);
\end{tikzpicture}
  \caption{The graphical model corresponding to the
    proposed construction.
    White circles represent unobserved hidden states
    and the black circle represents the observed
    data. Gray squares represent a jointly normal
    distribution. The arrow indicates a model interrogation.
    An implied non-Gaussian factor between $Y^{(0)}(\tn)$
    and $z_n$ is ignored to obtain a practical algorithm.}
  \label{fig:graphical-model}
\end{figure}
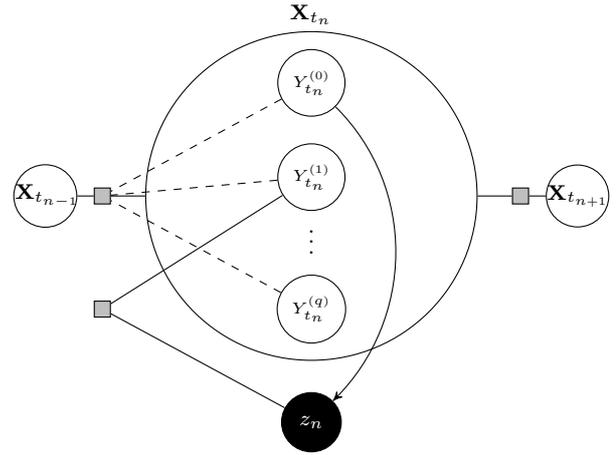

\editone{
Algorithm~\ref{alg:pseudo} conveys the general idea of a probabilistic ODE solver while
omitting parameter tuning aspects \del{that should be included in a modern solver} like
error control and step size selection.
The exact form of line~\ref{alg:line-evaluate}
depends on the choice of observation construction
and data likelihood model.
Without
data, the prior induces a probability distribution on
the hidden state $\vec{\edit{X}}_\tn$.
It remains to construct an observation $z_n$ and
a likelihood model $P(z_n\g\vec{X}_\tn)$.
}


\subsection{Observation assumptions}
\label{sec:observation-models}

\edit{
Recall from Sec.~\ref{sec:state-space} the
prior state space assumption
\begin{equation}
  \vec{X}_t = (Y_t, Y'_t, \dotsc, Y^{(q)}_t)\Trans \sim \N(\vec{m}_t, \mat{C}_t).
  \label{eq:state}
\end{equation}
Combining Equations~\eqref{eq:ivp-problem} and \eqref{eq:state} gives
  \begin{align}
    P(Y'_t) &= f(t, \cdot) \circ \N(Y_t; (\vec{m}_t)_0, (\mat{C}_t)_{00})\label{eq:push-forward}\\
    &\neq \N(Y'_t; (\vec{m}_t)_1, (\mat{C}_t)_{11})\label{eq:approximation}
  \end{align}
where Equation~\eqref{eq:push-forward} denotes the transformed random variable.
The exact form of that push-forward is not usually tractable
for general $f$ (one exception are linear ODEs, which of
course do not require nontrivial numerical algorithms).
}

\edit{
We will show below, however, that replacing the push-forward with an approximate
inference step captured by a Gaussian likelihood 
leads to good analytic properties
of the resulting Gaussian posterior.
This likelihood, which ignores the recursive nature of the ODE
(Equation~\eqref{eq:push-forward} and Figure~\ref{fig:graphical-model}), will be parametrized as
\begin{equation}
  P(z_n\g Y'_\tn) = \N(z_n; Y'_\tn, R_n^2)
  \label{eq:obs-model}
\end{equation}
where $z_n$ are the observations that have
yet to be constructed and $R_n^2$ can be interpreted
as an observation uncertainty.
Another way to phrase Equation~\eqref{eq:obs-model} is to write
\begin{equation}
  z_n = \mat{H}_1 \vec{X}_\tn + \nu
  \label{eq:obs-noise}
\end{equation}
where the latent variable $\nu\ce y'(\tn) - f(\tn,\mat{H}_0\vec{X}_\tn)$
captures the error between $f$ at the estimated solution
and the true solution's derivative.
 The approximation in Equation~\eqref{eq:obs-model} is to assign
 a centered Gaussian density $P(\nu)=\N(\nu; 0,R_n^2)$ to this latent variable. 
Purely from a formal perspective this
$\nu$ is a ``random variable'', but we stress again
that $P(\nu)$ captures \emph{uncertainty} arising
from lack of computational information
about a deterministic quantity, not any physical
sort of randomness in a frequentist sense.
That is, solving the same IVP several times will
of course always produce the exact same $\nu$,
but that same $\nu$ will always be just as unknown.
Figure~\ref{fig:graphical-model} displays a graphical
model corresponding to the construction.
All current probabilistic numerical ODE solvers share this particular assumption~\eqref{eq:obs-model}
\cite{skilling1991bayesian,o.13:_bayes_uncer_quant_differ_equat,%
schober2014nips,conrad_probability_2015,KerstingHennigUAI2016,teymur2016probabilistic}.
The differences between these algorithms chiefly
lies in the prior on $\vec{X}_t$,\
and how the observation $z_n$ is produced within the algorithm.
}

\edit{
It remains to construct $z_n$ and $R_n^2$.
One possible definition of $z_n$
}%
\editone{
is to evaluate
\begin{equation}
  \begin{aligned}
	z_n &\edit{\leftarrow} \Exp[f]\\
      &= \int f(\tn,Y_\tn)\,\N(Y_\tn; (\edit{\vec{m}}^-_\tn)_0, (\edit{\mat{C}}^-_\tn)_{00})\, \diff Y_\tn,
  \end{aligned}
	\label{eq:data-generation}
\end{equation}
where $\N(\vec{X}_\tn; \vec{m}^-_\tn, \mat{C}^-_\tn) = P(\vec{X}_\tn\g z_{[n-1]})$
is the prediction distribution of $\vec{X}_\tn$
given the data $z_{[n-1]}$
}%
\edit{
and $\leftarrow$ denotes assignment in code.  
}

\editone{
With these conventions, two new issues
emerge: the evaluation of the intractable
Equation~\eqref{eq:data-generation} and the
determination of $\edit{R_n^2}$. Kersting
and Hennig~\cite{KerstingHennigUAI2016} propose
to put
\begin{equation}
 \edit{R_n^2} \edit{\leftarrow} \int f(\tn,Y_\tn)^2 \,\N(Y_\tn; (\edit{\vec{m}}^-_\tn)_0, (\edit{\mat{C}}^-_\tn)_{00})\,\diff Y_\tn - \Exp[f]^2
 \label{eq:data-noise}
\end{equation}
and to evaluate both integrals by Bayesian
quadrature.
Chkrebtii et al.'s~\cite{o.13:_bayes_uncer_quant_differ_equat}
method draws a sample
$u_n \sim \N((\edit{\vec{m}}^-_\tn)_0, (\edit{\mat{C}}^-_\tn)_{00})$,
computes $z_n \edit{\leftarrow} f(\tn,u_n)$ and \edit{$R_n^2$}
is set to $(\edit{\mat{C}}^-_\tn)_{11}$.
In light of Kersting \& Hennig~\cite{KerstingHennigUAI2016},
this could be thought of as a form of
Monte Carlo scheme to evaluate \eqref{eq:data-generation}.
}

\edit{
As a further restriction to the likelihood~\eqref{eq:obs-model}
more widely used by other probabilistic numerical solvers,
we will here focus on models with $R_n^2\to 0$. That is
\begin{equation}
\begin{aligned}
  z_n &\leftarrow f(t_n, (\vec{m}^-_\tn)_0),\\
  P(z_n\g Y'_\tn) &= \delta(z_n - Y'_\tn) = \N(z_n; Y'_\tn, 0),
\end{aligned}
\label{eq:data-choice}
\end{equation}
This means the estimation node $y^-_\tn$ for the evaluation of $f$
is simply the current mean prediction, and the resulting observation is modeled as being correct.
}

\edit{
From the analytical viewpoint external to the algorithm itself,
of course, one does not expect that the model assumption
of a Gaussian likelihood, much less one with vanishing width, to hold in reality.
The point of the analysis in Sec.~\ref{sec:nordsieck} is to
demonstrate that this model and evaluation scheme yields 
a method of nontrivial convergence order for some choices
of state spaces, while simultaneously keeping computational
cost very low (that is very similar to that of classic multistep solvers).
}%
\editone{
That is because the predictive
posterior distributions $P(\vec{X}_\tn\g z_{[n]})$ can
be computed by the linear-time algorithm known as
\emph{Kalman filtering} \cite{kalman1960new,sarkka2006thesis,sarkka2013bayesian}.
The marginal predictive posterior distributions
given all data \edit{$P(\vec{X}_t\g z_{[N]})$} can be computed using the
\emph{Rauch--Tung--Striebel smoothing} equations
\cite{rauch1965maximum,sarkka2006thesis,sarkka2013bayesian}.
Simultaneously, one can draw samples from the
full joint posterior. These two operations
increase the computational cost minimally:
They require additional computations comparable to those
used for interpolation in classic solvers;
but neither smoothing nor sampling requires
additional evaluations of $f$.
The computational complexity stays linear in number of data points collected.
If the full joint posterior is also required for some reason,
this is also possible to construct \cite{solinPhDthesis,sdeToGP}.
}%
\edit{
As a second consequence, the computation
becomes deterministic which enables unit testing of the
resulting code.
}


\edit{
As a side remark, we note some obvious restrictions
of the combination of Gaussian (process) prior and
likelihood used here:
Since this combination means
the posterior is always a Gaussian process,
one cannot hope to accurately capture bifurcation events,
higher order correlations in the discretization errors,
or other higher order effects.
}






\subsection{Detailed example}
\label{sec:example}

Consider a concrete example.
We solve the following IVP
\begin{equation}
\begin{aligned}
  y' &= f(t,y) = f(y) = ry(1 - \nicefrac{y}{K}),\\
  y(t_0) &= y_0 = \nicefrac{1}{10},\q r = 3, K = 1,  
\end{aligned}
\label{eq:example}
\end{equation}
\editone{on the interval $[0, 1.5]$}.
Equation~\eqref{eq:example} is the sigmoid logistic growth function. Its solution
is available in closed form
\begin{equation*}
y(t) = \frac{K y_0 \exp(rt)}{K + y_0(\exp(rt) - 1)}.  
\end{equation*}

To solve this system, we apply a $2$-times
integrated Wiener process.
\editone{For this example, we fix $h_n = h = 0.3$,
such that $t_n = t_0 + hn$ for all $n$.}
Usually, the initial values are chosen to be $\vec{m}_\tmo^- = \vec{0}$
and $\mat{C}_\tmo^- = \mat{Q}(\infty)$, which is the
so-called steady state for stationary processes \cite{hartikainen2010kalman}.
The latter \edit{does not exist in the case} for the integrated
Wiener process, since the IWP is not stationary.
However, as has been shown in~\cite{schober2014nips},
this can be done analytically,
collecting the first $q$ derivative observations
$z_0, \dotsc, z_{q-1}$ manually
\editone{in the interval $[t_0, t_1]$}
and then inserting them in
the analytic formulas,
\editone{
yielding the filtering distribution
$P(\vec{X}_{t_1}\g z_{[q-1]}) = \N(\edit{\vec{X}_{t_1};}\,\vec{m}_{t_1}, \mat{C}_{t_1})$
(see also Sec.~\ref{sec:connection-rk}).
The remainder of the interval $[t_1, \dotsc, t_N = T]$
is solved with the familiar Kalman filter equations
\begin{align}
\vec{m}_\tn^- &= \mat{A}(h) \vec{m}_\tnmo, \label{eq:pred-mean}\\
\mat{C}_\tn^- &= \mat{A}(h) \mat{C}_\tnmo \mat{A}(h)\Trans + \mat{Q}(h) \label{eq:pred-cov}\\
\intertext{and}
\lambda_n &= f(\edit{\tn}, \mat{H}_0 \vec{m}_\tn^-) - \mat{H}_1 \vec{m}_\tn^-, \label{eq:kalman-res}\\
\mat{K}_n &= \mat{C}_\tn^-\mat{H}_1\Trans[\mat{H}_1 \mat{C}_\tn^- \mat{H}_1\Trans]^{-1}, \label{eq:kalman-gain}\\
\vec{m}_\tn &= \vec{m}_\tn^- + \mat{K}_n \lambda_n, \label{eq:update-mean}\\
\mat{C}_\tn &= \mat{C}_\tn^- - \mat{K}_n [\mat{H}_1 \mat{C}_\tn^- \mat{H}_1\Trans]\mat{K}_n\Trans. \label{eq:update-cov}
\end{align}
}

\begin{figure*}
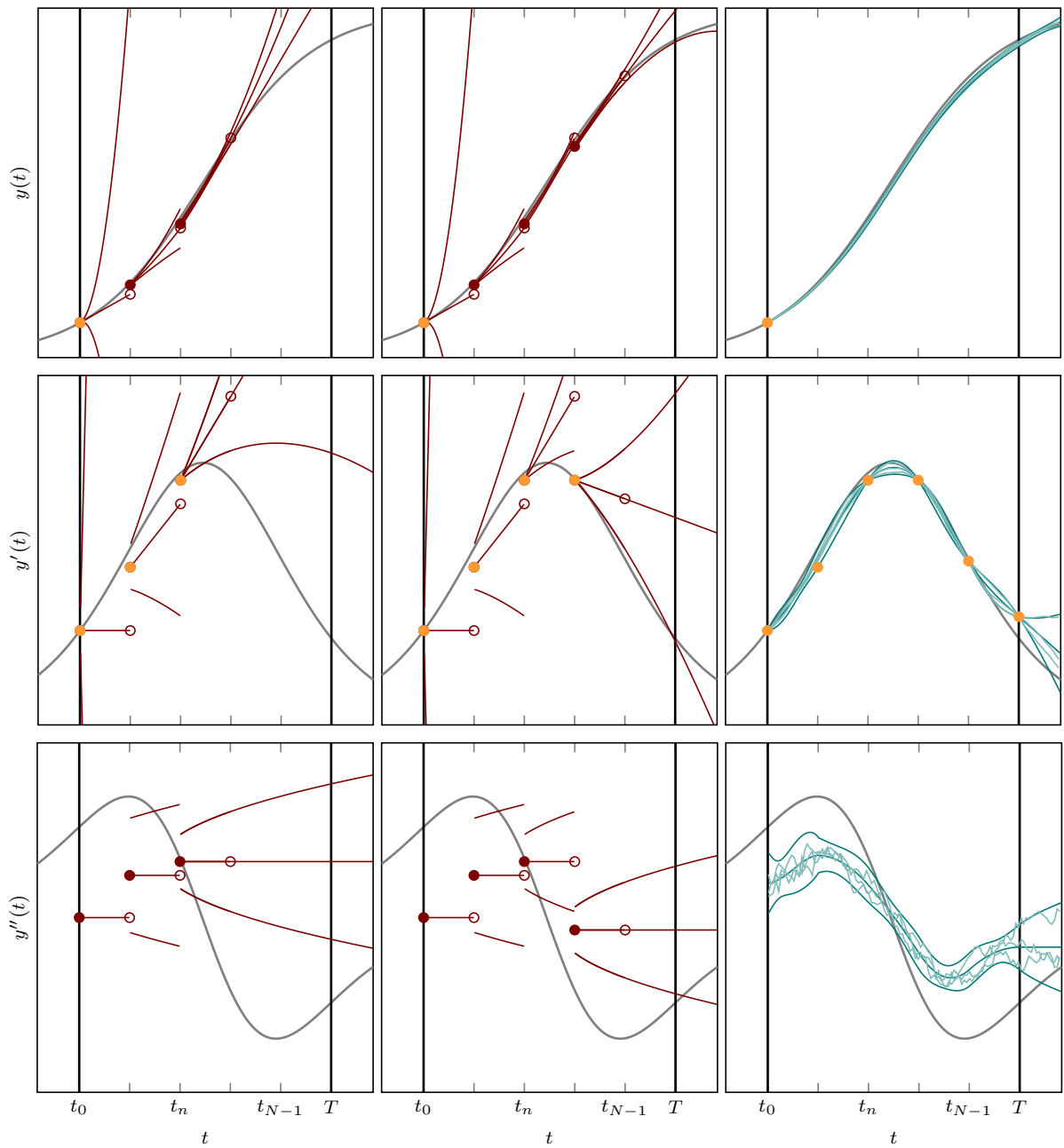

  \centering
  \footnotesize
  \setlength{\figurewidth}{0.3\textwidth}
  \setlength{\figureheight}{\figurewidth}
  \input{sec2_example-yt-4p}
  \input{sec2_example-yt-5p}
  \input{sec2_example-yt-smoothed}\\
  \input{sec2_example-dyt-4p}
  \input{sec2_example-dyt-5p}
  \input{sec2_example-dyt-smoothed}\\
  \input{sec2_example-ddyt-4p}
  \input{sec2_example-ddyt-5p}
  \input{sec2_example-ddyt-smoothed}
  \caption{The $2$-times integrated Wiener process
    $\diff \vec{X}_t = \mat{U}_3 \vec{X}_t \diff t + \vec{e}_2 \diff W_t$ applied to the
    logistic growth problem $y' = ry(1 - y/K)$. The plot shows the true solution
    (gray) of the function $y$ and its first two derivatives, as well as the
    numerical solution $Y$, given by its mean $\vec{m}_i$ (red line) and
    covariance $\mat{C}$, visualized as point-wise plus-minus twice the standard
    deviation $\vec{m}_i \pm 2 \sqrt{\mat{C}_{ii}}$. Empty circles are
    predicted values at time $t_n$, filled circles represent updated values, orange
    dots are function and derivative observations.
    The first two columns display two \emph{predict-evalute-update-predict} cycles.
    The last column shows the smoothed final solution (green, thick lines) and three
    samples from the predictive posterior (thin lines).}
  \label{fig:example-filter}
\end{figure*}

Figure~\ref{fig:example-filter} shows the state of
the algorithm after $2$ steps have been taken.
The solution looks discontinuous, because the
information of later updates $z_n$ has not been
propagated to previous time points $t_m, m < n$.
The last column of Figure~\ref{fig:example-filter} shows the
\emph{(predictive posterior) smoothing distribution}
wherein all the information is globally available.


\section{Classical analysis for the probabilistic method}
\label{sec:connection}

\editone{
The most important test for any numerical algorithm is
that it works in practice and under the requirements of
potential users.
The proposed probabilistic numerical algorithm has been
motivated and derived from the computational properties
that established classical algorithms provide.
The classical algorithms have}
been studied intensely for over a century, to a point where the
theory could almost be considered complete \cite{gear1981anything}.
Thus, a newly proposed algorithm---even when motivated from a 
different back\-ground---should stand up to classical
analysis.

\editone{
While many specialized models and algorithms have been
proposed, two}
standard classes of algorithms have become
prevalent for the solution of \eqref{eq:ivp-problem}: \emph{Runge--Kutta (RK)} methods
and \emph{(linear) multistep methods (LMMs)} or combinations thereof (general
linear methods, GLMs \cite{Butcher1985general}).
These classes share a similar type of algorithmic structure and analysis:
at time $\tn$, evaluate $f$ with a numerical approximation\del{s} $\edit{y_\tn}$ to construct
an updated numerical approximation $\edit{y_\tnpo}$ from linear combinations of
the function evaluations $\edit{f_\tn}$ \edit{(exact definitions below)}.
The update weights are parameters
of a given method and, if chosen appropriately, can be shown to coincide
with the Taylor approximation of the true solution $y$ up to $q$ terms.

\edit{In the following, we present} results relating the
newly proposed probabilistic method to existing
algorithms,
\editone{%
which allows us to transfer the known
results to our method.
Interpreting the probabilistic model
from the viewpoint of classical analysis
adds a justification to the assumptions
made in the previous sections
by saying that these assumptions---unintuitive
as they may be at first---are the same assumptions
that are implied by the application of
a classical algorithm.
}

\subsection{On the connection to Nordsieck methods}
\label{sec:nordsieck}
\editone{
Linear multistep methods are defined by the relationship
\begin{equation}
  \sum_{i=0}^q \alpha_i y_{t_{n-i}} = h \sum_{i=0}^q \beta_i f_{t_{n-i}},
  \label{eq:multistep}
\end{equation}
where $f_\tn$ are approximations to $y'(\tn)$, $h$ is the step size,
and the $\alpha_i$ and
$\beta_i$ are parameters of the method. If $\beta_0 = 0$,
then \eqref{eq:multistep} defines $y_\tn$ without requiring $f_\tn$
and we can set $f_\tn \ce f(\tn,y_\tn)$ for the computation of $y_\tnpo$.
This is called an \emph{explicit method}.
In contrast, if $\beta_0 \neq 0$, we still define $f_\tn = f(\tn,y_\tn)$,
but \eqref{eq:multistep} now defines a non-linear equation for
$y_\tn$ given a non-linear $f$. We say that $y_\tn$ is \emph{implicitly defined}
and, therefore, methods with $\beta_i \neq 0$ are called \emph{implicit methods}.
Assuming that at least one of $\{\alpha_q,\beta_q\}$ does not vanish,
the method requires the numerical approximation on $q$ previous grid
points and \eqref{eq:multistep} is called a $q$-step method.
}

\editone{Skeel showed in \cite{skeel1979equivalent} that implicit} LMMs
can be written in Nordsieck form \cite{nordsieck1962numerical}:
\begin{align}
\vec{x}_{n}   &= \edit{\left(y_\tn, hy'_\tn, \dotsc, \nicefrac{h^qy^{(q)}_\tn}{q!}\right)\Trans},\label{eq:nordsieck-vector}\\
\vec{x}_{n+1} &= \left(\Id - \vec{l}\vec{e}_1\Trans\right)\mat{P} \vec{x}_n + h \vec{l} z_n,\label{eq:nordsieck}
\end{align}
where $\mat{P}$ is the Pascal triangle matrix with entries 
$p_{ij} = \ind{i \le j}\binom{j}{i}$ and
$\vec{l} \editone{= \vec{l}(\{\alpha_i,\beta_i\})}$ is the weight vector defining the method.
\editone{
The vector $\vec{x}_n$ is called \edit{the} \emph{Nordsieck vector} in
honor of its inventor Nordsieck \cite{nordsieck1962numerical},
and a LMM written in Nordsieck form is also called a \emph{Nordsieck method}.
The intuition behind this family of methods is to replace $y(t)$
or $y'(t)$ on $[t_{n-q}, \dotsc, \tn]$ with a local polynomial
$\pi(t) = \pi_{[t_{n-q}, \dotsc, \tn]}(t)$ of order $q$.
}

\editone{
The difference in presentation between \eqref{eq:multistep}
and \eqref{eq:nordsieck-vector} can be understood as expressing
$\pi(t)$ either in Lagrange notation (Equation~\eqref{eq:multistep})
or Taylor expansion notation (Equation~\eqref{eq:nordsieck-vector}).
\edit{In this case}, $\mat{P}\vec{x}_n$ yields a prediction of the numerical
Taylor expansion at $\tnpo$ and the%
}
scalar increment $z_n$ is implicitly defined as the solution to
\begin{equation}
\editone{h^{-1}(\mat{P}\vec{x}_n)_1 + l_1 z_n = f(\tn + h, (\mat{P}\vec{x}_n)_0 + hl_0z_n)},\label{eq:implicit-method}
\end{equation}
\editone{which is the correction from $\vec{x}_n$ to $\vec{x}_{n+1}$
to the Taylor coefficients.
Equation~\eqref{eq:implicit-method} can be solved by iterated function
evaluations of the form
\begin{align}
  z_n^{(1)} &\ce f(t_n + h, (\mat{P}\vec{x}_n)_0)\label{eq:pec-1}\\ 
  z_n^{(M)} &\ce l_1^{-1} [f(t_n +h, (\mat{P}\vec{x}_n)_0 + hl_0z_n^{(M-1)})\notag\\
              &\phantom{\ce l_1^{-1}[f]}
              - h^{-1}(\mat{P}\vec{x}_n)_1]\label{eq:pec-inf}
\end{align}
or by directly solving \eqref{eq:implicit-method} with some
variant of the Newton-Raphson method.
}

If \editone{$z^{(M)}$ is used in the computation of \eqref{eq:nordsieck},
the resulting algorithm} is called a P(EC)\textsuperscript{M} method.
If Equation~\eqref{eq:implicit-method} is
solved up to numerical precision, the method is called a P(EC)\textsuperscript{$\infty$} method.
Nordsieck methods with suitable weights $\vec{l}$
can be shown to have local truncation error
of order $q$ or $q+1$
\cite{skeel1977consistency,skeel1979equivalent}.
More details can also be found in standard textbooks
\cite{hairer87:_solvin_ordin_differ_equat_i,DeuflhardScientific}.
\del{These methods are the only ones amenable
to the description offered below, but they will serve as sufficiently
interesting connections to the statistical construction.}

\editone{
We will now show how the Kalman filter \eqref{eq:pred-mean}--\eqref{eq:update-cov}
can be rewritten such that the mean prediction takes the form of
\eqref{eq:nordsieck}. This enables to analyze the
proposed algorithm in light of classical Nordsieck method results,
but can also guide the further development of the probabilistic approach
with the experience of existing software.
}

Considering a fixed step size $h_n = h,\,n = 1, \dotsc, N$,
we re-scale the state space and SDE of the IWP($q$) by scaling matrix
$\mat{B}$ to define an equivalent notation
\begin{equation}
\begin{aligned}
  \vec{\tilde{X}}\edit{_t} &= 
  \begin{pmatrix}Y_t, & h Y'_t, & \frac{h^2}{2!} Y''_t, & \hdots & \frac{h^q}{q!} Y^{(q)}_t\end{pmatrix}\Trans\\
  &= \underbrace{\operatorname{diag}\begin{pmatrix}1, & h, & \frac{h^2}{2!}, & \hdots & \frac{h^q}{q!}\end{pmatrix}}_{\ec \mat{B}} \vec{X}_t,\label{eq:rescaled-state}
\end{aligned}.
\end{equation}
This state vector is the \emph{Nordsieck vector}. The advantage of this notation is that
\eqref{eq:sde} simplifies to
\begin{equation}
  \diff \vec{\tilde{X}}_{\edit{t}} = \mat{B} \mat{U}_{q+1} \mat{B}^{-1} \vec{\tilde{X}}_{\edit{t}} \diff t+ \mat{B} \vec{e}_q \diff W_{\edit{t}},
  \label{eq:sde-tilde}
\end{equation}
where $\mat{\tilde{A}}(h) = \mat{P}$, the Pascal triangle matrix,
and 
\begin{align}
(\mat{\tilde{Q}}(h))_{ij} &= (\mat{B} \mat{Q}(h) \mat{B}\Trans)_{ij}\notag\\
&= \frac{h^i}{i!} \ssq  \frac{h^{2q+1-i-j}}{(2q+1-i-j)(q-i)!(q-j)!} \frac{h^j}{j!}\notag\\
&= \frac{\ssq h^{2q+1}}{(2q+1-i-j)(q-i)!(q-j)!i!j!}
\label{eq:iwp-tilde-diffusion}
\end{align}
which can be seen by inserting \eqref{eq:sde-tilde} into
\eqref{eq:sde-sol-mu-cov} and simplifying.
Furthermore, the observation matrices become
$\mat{\tilde{H}}_0 = \mat{H}_0 \mat{B}^{-1} = \vec{e}_0$
and $\mat{\tilde{H}}_1 = \mat{H}_1 \mat{B}^{-1} = h^{-1} \vec{e}_1$.
Rewriting the filtering equations, we arrive at
\begin{align}
\mat{C}_\tn^- &= \mat{P} \mat{C}_\tnmo \mat{P}\Trans + \mat{\tilde{Q}}(h),\label{eq:pred-cov-nordsieck}\\
\mat{K}_n &= \mat{C}_\tn^-\mat{\tilde{H}}_1\Trans[\mat{\tilde{H}}_1 \mat{C}_\tn^- \mat{\tilde{H}}_1\Trans]^{-1}\label{eq:kalman-gain-nordsieck}\\
\intertext{and}
\vec{m}_\tn &= (\Id - \mat{K}_n \mat{\tilde{H}}_1) \mat{P}\,\vec{m}_\tnmo + \mat{K}_n \editone{z_n^{(1)}},\label{eq:update-mean-nordsieck}\\
\mat{C}_\tn &= (\Id - \mat{K}_n \mat{\tilde{H}}_1) \mat{P} (\mat{C}_\tnmo \mat{P}\Trans + \mat{\tilde{\Phi}}_{12}(h)\Trans)\label{eq:update-cov-nordsieck}.
\end{align}
Choosing a prior covariance matrix \editone{with entries $(\mat{C}_\tmo^-)_{ij} = \ssq h^{2q+1}c_{ij}$},
for some $c_{ij} \in \Re$ such that \editone{$\mat{C}_\tmo^-$} is a valid covariance matrix,
it can be shown by induction that all entries of $\mat{C}_\tn$ for all $n$ have this
structural form.
As a by-product, $\mat{K}_n = h\,(k_{n,0}, 1, k_{n,2}, \dotsc, k_{n,q})\Trans$ for some
$k_{n,i} \in \Re$ which follows from \eqref{eq:kalman-gain-nordsieck}.

Given these invariants, Equation~\eqref{eq:update-mean-nordsieck}
has the structure of a multistep method written
in Nordsieck form \eqref{eq:nordsieck}.
\editone{
The only difference is the changing weight vector $\vec{K}_n$~%
\eqref{eq:update-mean-nordsieck} as compared
to the constant weights in \eqref{eq:nordsieck}.
Multistep methods with varying weights have been
studied in the literature \cite{crouzeix1984convergence,brown1989vode}.
\edit{These works are often in the context of
variable step sizes $h_n \neq h$,
but variable-coefficient methods have also been
studied for other purposes, for example, cyclic methods \cite{Albrecht1978}.
These works have in common that the weights are
free variables that are not limited through the
choice of model class.
As a consequence, determining optimal weights
can be algebraically difficult \cite[\textsection III.5]{hairer87:_solvin_ordin_differ_equat_i}.
}
}

\editone{
Here, variable step sizes are easily obtained
by working with the representation \eqref{eq:sde}
instead of \eqref{eq:sde-tilde} and
computing \eqref{eq:trans-and-diffusion}
according to $h_n$.
\edit{
In contrast to classical methods,
the weights $\vec{K}_n$ cannot be chosen freely,
but are determined through the choice of model~\eqref{eq:sde},
and the evolution of
}%
the underlying uncertainty $\mat{C}_\tn$.
While Kersting and Hennig~\cite{KerstingHennigUAI2016}
provide some preliminary empirical evidence that
these adaptive weights $\vec{K}_n$ might
actually improve the estimate, more
rigorous analysis is required for theoretical
guarantees.
}

\editone{
In fact, Skeel (and Jackson) \cite{skeel1976analysis,skeel1977consistency}
consider more general \emph{propagation matrices} $\mat{S}$ for
$\edit{\vec{x}_\tn} = \mat{S} \edit{\vec{x}_\tnmo}$ in Equation~\eqref{eq:nordsieck}.
\edit{
Every model of the form~\eqref{eq:sde} implies such a
general propagation matrix by identifying
$\mat{S}_n = (\Id - \vec{K}_n \mat{H}_1) \mat{A}(h_n)$.
Thus, applying the Kalman filter to LTI SDE models
is structurally equivalent to a variable-coefficient
multistep method. 
}%
This motivates the following%
}
definition
and Algorithm~\ref{alg:prob-nordsieck} for the
probabilistic solution of initial value problems.

\begin{definition}
\label{def:prob-nordsieck}
A \editone{\emph{probabilistic filtering ODE solver (PFOS)}} is the Kalman
filter applied to an initial value problem with an underlying
Gauss--Markov linear, time-invariant SDE \del{model }and Gaussian observation
likelihood \editone{model}.
\end{definition}

\newcommand{\widevar}[1]{\mathrlap{#1}\phantom{\vec{m}_\tmo}}

\begin{algorithm}
\caption{Probabilistic \editone{filtering ODE solver}}
\label{alg:prob-nordsieck}
\begin{algorithmic}[1]
\STATE{Define $\tmo \ce t_0$, choose $\mat{F}, \vec{L}\edit{, \ssq}$, initialize $\vec{m}_\tmo^-, \mat{C}_\tmo^-$ accordingly}
\STATE{$\widevar{\vec{K}_{-1}} \gets \mat{C}_\tmo^-\mat{H}_0\Trans [\mat{H}_0\mat{C}_\tmo^-\mat{H}_0\Trans]^{-1}$}
\STATE{$\widevar{\vec{m}_\tmo} \gets \vec{m}_\tmo^- + \vec{K}_{-1} [y_0 - \mat{H}_0 \vec{m}_\tmo^-]$}
\STATE{$\widevar{\mat{C}_\tmo} \gets \mat{C}_\tmo^- - \vec{K}_{-1} [\mat{H}_0\mat{C}_\tmo^-\mat{H}_0\Trans] \vec{K}_{-1}\Trans$}
\FOR{n = 0 \TO N}
\STATE{$\widevar{h_n} \gets \tn - \tnmo$}
\STATE{Compute $\mat{A}(h_n),\,\mat{Q}(h_n)$}
\STATE{$\widevar{\vec{m}_\tn^-} \gets \mat{A}(h_n) \vec{m}_\tnmo$}\hfill\COMMENT{Predict}
\STATE{$\widevar{\mat{C}_\tn^-} \gets \mat{A}(h_n) \mat{C}_\tnmo \mat{A}(h_n)\Trans + \mat{Q}(h_n)$}
\STATE{$\widevar{z_n} \gets f(t_n, \mat{H}_0 \vec{m}_\tn^-)$}\hfill\COMMENT{Evaluate}\label{alg:psof-eval}
\STATE{$\widevar{\lambda_n} \gets z_n - \mat{H}_1 \vec{m}_\tn^-$}\hfill\COMMENT{Update}
\STATE{$\widevar{\vec{K}_n} \gets \mat{C}_\tn^- \mat{H}_1\Trans [\mat{H}_1\mat{C}_\tn^-\mat{H}_1\Trans]^{-1}$}\label{alg:psof-up}
\STATE{$\widevar{\vec{m}_\tn} \gets \vec{m}_\tn^- + \vec{K}_n \lambda_n$}
\STATE{$\widevar{\mat{C}_\tn} \gets \mat{C}_\tn^- - \vec{K}_n [\mat{H}_1\mat{C}_\tn^-\mat{H}_1\Trans]^{-1}\vec{K}_n\Trans$}
\ENDFOR
\RETURN \edit{$\{\vec{m}_{t_N}, \mat{P}_{t_N},\; n = -1,\dotsc,N \}$}
\end{algorithmic}
\end{algorithm}

\editone{
As was the case in Algorithm~\ref{alg:pseudo}, the exact form of
lines~\ref{alg:psof-eval}--\ref{alg:psof-up} depend on the
choice of likelihood model (cf. Kersting and Hennig \cite{KerstingHennigUAI2016}).
}

\editone{
We will now study the long-term behavior of the PFOS. In particular,
we will ask what is the long-term behavior for the sequence of
Kalman gains $(\mat{K}_n)_{n=0,\dotsc}$ and how this will influence
the solution quality.
It can be shown that its
properties are linked to properties of the
\emph{discrete algebraic Riccati equation}, of
which the theory has largely been developed \cite{lancaster1995algebraic}.
\edit{
Denote by
$\mat{\gamma} : \Re^{(q+1) \times (q+1)} \to \Re^{(q+1) \times (q+1)}$
the function that maps the covariance matrix $\mat{C}_\tnmo$
of the previous knot $\tnmo$ to the covariance matrix $\mat{C}_\tn$
at the current knot $\tn$ (Equation~\eqref{eq:update-cov-nordsieck}).
}
If there exists a (unique) fixed point $\mat{C}^*$
of $\mat{\gamma}$, it's called the \emph{steady state}
of the model~\eqref{eq:sde}.
Associated with a fixed point $\mat{C}^*$ is also
a constant Kalman gain $\mat{K}^*$ that is obtained
at the (numerical) convergence of $\mat{C}^*$.
}

We will now show that
\editone{there is a subset of model~\eqref{eq:sde} that
converges to a steady state. This subsystem completely
determines a constant Kalman gain $\mat{K}^*$}
at least in the case of
the IWP($1$) and IWP($2$). Thus, just like in the equivalence result
for the Runge--Kutta methods
\editone{in Schober et al.~\cite{schober2014nips},
the \edit{result} of the PFOS} is equivalent (in the sense of numerically identical)
\editone{after an initialization period} to
a corresponding classical Nordsieck method defined by the weight
vector $\editone{\mat{K}^*}$ and we can apply all the known theory of
multistep methods to the mean of the \editone{PFOS}.

\begin{proposition}
The PFOS arising from the
once integrated Wiener process IWP$(1)$ is equivalent in
its predictive posterior mean 
\editone{to the P(EC)\textsuperscript{1} implementation of}
the trapezoidal rule.
\label{prop}
\end{proposition}
\editone{
\begin{proof}
The trapezoidal rule, written as in Equation~\eqref{eq:multistep}, is
\begin{equation}
  y_\tn = y_\tn + \frac{h}{2} (f_\tnmo + f_\tn).\label{eq:trap-rule}
\end{equation}
We will show that $(\vec{m}_\tn)_0 = (\vec{m}_\tnmo)_0 + \nicefrac{h}{2}[(\vec{m}_\tnmo)_1 + (\vec{m}_\tn)_1)]$
for all $n$ by induction.
Let $\vec{m}_\tmo^- = \vec{0}$ and
$\mat{C}_\tmo^- \in \Re^{2\times 2}$ be an arbitrary
covariance matrix. Applying the first three lines of Algorithm~\ref{alg:prob-nordsieck}
algebraically, the predicted values are
\begin{align}
\vec{m}_\tmo = &\begin{pmatrix} y_0 \\ m_{t_0,1}^-\end{pmatrix},
& \mat{C}_\tmo = &\begin{pmatrix} 0 & 0 \\ 0 & c_{t_0, 11}^-\end{pmatrix}
\end{align}
for some $m_{t_0,1}^-, c_{t_0, 11}^-$.
Continuing in this fashion yields
$z_0 \ce f(t_0, y_0)$ and $\vec{m}_{t_0} = (y_0, z_0)\Trans, \mat{C}_{t_0} = \mat{0}$.
Using \eqref{eq:pred-mean} and \eqref{eq:pred-cov} to compute the
predictions at $t_1$, we arrive at
\begin{align}
  \vec{m}_{t_1}^- &= \begin{pmatrix} y_0 + h z_0 \\ z_0\end{pmatrix},
  & \mat{C}_{t_1}^- &= \mat{Q}(h)
  \label{eq:explicit-euler}
\end{align}
and we see that $\mat{H}_0 \vec{m}_{t_0 + h}^- = y_0 + h z_0 = (\mat{P}(y_0, h z_0)\Trans)_0$.
Completing the Kalman step by applying
Equations~\eqref{eq:kalman-res}-\eqref{eq:update-cov} yields 
\begin{align}
\vec{m}_{t_1} &= \begin{pmatrix} y_0 + \frac{h}{2}[z_0 + z_1]\\ z_1\end{pmatrix},
& \mat{C}_{t_1} &= \ssq \begin{pmatrix} \frac{h^3}{12} & 0 \\ 0 & 0\end{pmatrix},
\label{eq:heuns-method}
\end{align}
where $z_1 \ce f(t_1, y_0 + h z_0)$.
Comparing with \eqref{eq:pec-1}, we see that $z_1$ is of the
desired form, which completes the start of the induction.
Finally, we observe that the second column of
$\mat{C}_{t_1} = \vec{0} = \mat{C}_{t_0}$, i.e, this will be invariant
and, as a consequence, the second column of $\mat{C}^-_\tn$ is simply the
second column of $\mat{Q}(h)$, and the induction is completed.
\smartqed\qed
\end{proof}
}
\del{This combination of results is in agreement with Skeel~\cite{skeel1979equivalent},
who notes that the iterative solution of the implicit
system \eqref{eq:update-mean-nordsieck} in Nordsieck form is
generally not equivalent to any P(EC)\textsuperscript{M}E
method which is the RK interpretation of the explicit trapezoidal rule.}


We will now investigate the model corresponding to the IWP($2$).
To this end, let $\mat{C}_\tnmo$ be of the form
\begin{equation}
   \mat{C}_\tnmo = \ssq h^5 \begin{pmatrix} c_{00} & 0 & c_{02}\\ 0 & 0 & 0\\ c_{02} & 0 & c_{22}\end{pmatrix}
   \label{eq:iwp-2-prior-cov}
\end{equation}
such that $\mat{C}_\tnmo$ is a valid covariance matrix.
Equations~\eqref{eq:kalman-gain-nordsieck} and \eqref{eq:update-cov-nordsieck}
guarantee that this form holds for $t_0$, that is, after
the first \edit{derivative} observation and this can be found
through the inspection of $\mat{\tilde{Q}}(h)$.
Applying one Kalman step algebraically, by inserting
\eqref{eq:iwp-2-prior-cov} into
\eqref{eq:kalman-gain-nordsieck} and \eqref{eq:update-cov-nordsieck},
we find that
\begin{equation}
\begin{aligned}
(\mat{C}_\tn)_{00} &= \ssq h^5 \tfrac{3840 c_{00} c_{22} + 320 c_{00} - 3840 c_{02}^2 + 110 c_{02} + 32 c_{22} + 1}{320 (12 c_{22} + 1)}\\
(\mat{C}_\tn)_{02} &= \ssq h^5 \tfrac{-(48 c_{02} + 24 c_{22} + 1)}{96 (12 c_{22} + 1)} = (\mat{C}_\tn)_{20}\\
(\mat{C}_\tn)_{22} &= \ssq h^5 \tfrac{16 c_{22} + 1}{16 (12 c_{22} + 1)}\\
(\mat{C}_\tn)_{ij} &= 0,\q i,j = 0, 1, 2, i \vee j = 1
\end{aligned}
\label{eq:iwp-2-post-cov}
\end{equation}
which has the same structural form.
\edit{
We will now consider the behavior of the coefficients $c_{ij}$.
Consider the dynamical system
$\bar{\gamma}_{22}(c) = (16c + 1)[16(12c+1)]^{-1}$
which maps the coefficient of the last entry in $\mat{C}_\tnmo$
to the next. The range and image of $\bar{\gamma}_{22}$
are the non-negative reals, since variances cannot be negative.
On this domain, $\bar{\gamma}_{22}$ has a continuous and bounded
derivative $\abs{\bar{\gamma}'_{22}} \leq \frac{1}{4}$.
In particular, $\bar{\gamma}_{22}$ is a contraction with
Lipschitz constant $\frac{1}{4}$.
Thus, the entries converge to the fixpoint $c^*_{22} = \frac{\sqrt{3}}{24}$
(which can be found with some simple algebra).
Similarly, one can either insert $c^*_{22}$ into
the respective form of $\bar{\gamma}_{02}$ or one
considers the two-dimensional mapping of both entries.
In both cases, a similar argument guarantees the
convergence to a fixpoint, which is found to be $c_{02}^* = - \frac{\sqrt{3}}{144}$.
}
Inserting these into Equation~\eqref{eq:kalman-gain-nordsieck},
we find that $\mat{K}_n = \mat{K}^* = (\frac{3 + \sqrt{3}}{12}, 1, \frac{3 - \sqrt{3}}{2})\Trans$
is the static probabilistic Nordsieck method of the IWP$(2)$ filter.
Inserting these weights into \cite[Theorem 4.2]{skeel1979equivalent}
\editone{proves the following:}
\begin{theorem}
The predictive posterior mean of the IWP($2$) with fixed
step size $h$ is a third order Nordsieck method\edit{, when
the predictive distribution has reached the steady state}.
\label{thm}
\end{theorem}

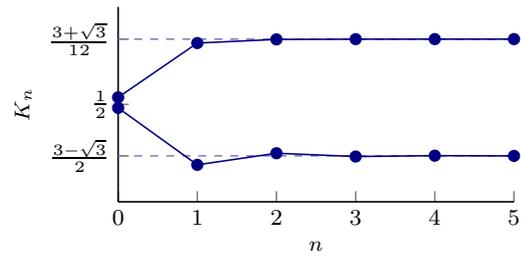
\begin{figure}
  \centering
  \setlength{\figurewidth}{0.45\textwidth}
  \setlength{\figureheight}{0.33\figurewidth}
%
%
\definecolor{mycolor1}{rgb}{0.50000,0.50000,0.75450}%
\definecolor{mycolor2}{rgb}{0.00000,0.00000,0.50900}%
\begin{tikzpicture}

\begin{axis}[%
width=0.665\figurewidth,
height=\figureheight,
at={(0\figurewidth,0\figureheight)},
scale only axis,
xmin=0,
xmax=5,
xlabel={$n$},
ymin=0.3,
ymax=0.7,
ytick={0.394337567297406,0.5,0.633974596215561},
yticklabels={{$\frac{3-\sqrt{3}}{2}$},{$\frac{1}{2}$},{$\frac{3+\sqrt{3}}{12}$}},
ylabel={$K_n$},
axis background/.style={fill=white},
axis x line*=bottom,
axis y line*=left,
mystyle
]
\addplot [color=mycolor1,dashed,forget plot]
  table[row sep=crcr]{%
0	0.394337567297406\\
5	0.394337567297406\\
};
\addplot [color=mycolor1,dashed,forget plot]
  table[row sep=crcr]{%
0	0.633974596215561\\
5	0.633974596215561\\
};
\addplot [color=mycolor2,solid,mark=*,mark options={solid},forget plot]
  table[row sep=crcr]{%
0	0.492808314288538\\
1	0.375905472588207\\
2	0.399710108777878\\
3	0.392931097251751\\
4	0.394716765793287\\
5	0.39423612984461\\
};
\addplot [color=mycolor2,solid,mark=*,mark options={solid},forget plot]
  table[row sep=crcr]{%
0	0.514383371422924\\
1	0.625905472588208\\
2	0.633397753604916\\
3	0.633933193578234\\
4	0.633971623705896\\
5	0.633974382799309\\
};
\end{axis}
\end{tikzpicture}%
  \caption{The weights $(K_n)_0$ and $(K_n)_2$ for $n = 0, \dotsc, 5$. 
  }
  \label{fig:weights}
\end{figure}

\edit{
Although Theorem~\ref{thm} is only valid
when the system has reached its steady state,
we find that the convergence (visualized in Figure~\ref{fig:weights}) is rapid in practice.
In the extreme case of $q=1$ (not shown),
in fact it is instantaneous, and Proposition~\ref{prop} is valid from the second step onwards.
This limitation could also be circumvented in practice
by initializing $\mat{C}_\tmo$ at steady state coefficients,
but this possibility is not required to achieve high-order
convergence on the benchmark problems we considered.}

\edit{
Figure~\ref{fig:weights} shows the situation for a
constant value of the diffusion amplitude $\sigma^2$.
In Section~\ref{sec:error-estimation}, we will discuss error estimation
and step size adaptation.
This process leads to a continuous adaptation of this variable,
which in turn means that the convergence shown in the figure
continues throughout the run of the algorithm.
So the practical algorithm presented here and empirically
evaluated in Section~\ref{sec:experiments} is not formally
identical to Nordsieck methods, merely conceptually closely related.
}

\editone{Inspecting the weights}
of the IWP($2$)\del{ to the Adams-Moulton
method of the same order}, we find that
\editone{this method has not been considered previously in the literature,
and, in particular, cannot be related to any of the typical formulas,
such as Adams-Moulton or backward differentiation formulas.}
This is not surprising, since the
\editone{\edit{result} of this method has been constructed to be twice
continuously differentiable},
whereas there is no such
guarantee for the solution provided by \editone{the typical methods}. In fact,
the probabilistic Nordsieck method is much closer related to
spline-based multistep methods such as
\cite{Loscalzo1967Spline,loscalzo1969introduction,Byrne1972Linear,andria1973integration}
since Gaussian process regression models have a one-to-one
correspondence to spline smoothing in a reproducing kernel Hilbert space
of appropriate choice \cite{kimeldorf1970correspondence,wahba1990spline}.
This also justifies the application of a full-support distribution,
even though it is known that the solution will remain in a compact set.
In the former case, the interpretation is one of average case error
whereas in the latter, the bound corresponds to the worst case error
\cite{paskov1993average}.

\edit{
The forms of $C_\tn$ found in Equations~\eqref{eq:heuns-method} \& \eqref{eq:iwp-2-post-cov}
also show that the standard deviation
$\operatorname{std}[(\vec{m}_\tn)_0] = (\mat{C}_\tn)_{00}^{1/2}$ can be meaningfully,
if weakly, interpreted as an approximation 
to the \emph{expected error} $\abs{y_\tn - y(\tn)}$ of the numerical solution
in the following local, asymptotic sense:
From our analysis of the IWP($q$), $q \in \{1, 2\}$,
we have $\abs{y_\tn - y(\tn)} \leq Ch^{q+1}$
whereas $(\mat{C}_\tn)_{00}^{1/2} \in \mathcal{O}(\sigma h^{q+1/2})$.
Estimating the intensity $\sigma$ of the stochastic process
amounts to estimating the unknown constant $C$.
}


We conclude this section by
\editone{considering some implications of the probabilistic interpretation
in contrast to}
other classical multistep methods.

\editone{Keeping all hyper-parameters (order $q$, prior diffusion
intensity $\ssq$, and step size $h$) fixed, the gain
$\mat{K}_n$ is completely determined, and, as a consequence,
we could have chosen to fully solve the implicit
equation~\eqref{eq:implicit-method} for the generation of $z_n$.
Solving \eqref{eq:implicit-method} up to numerical precision
can be interpreted to learn the true value
of the model \eqref{eq:sde} at $t_n$ which gives another
justification for using $R_n^2 = R^2 = 0$.
Since the P(EC)\textsuperscript{$\infty$} and
the P(EC)\textsuperscript{M} have the same order for all $M$
\cite{DeuflhardScientific}, this argument can
be extended to the case of the PEC\textsuperscript{1}
implementation which gives the most natural
connection to the Kalman filter.
}

\editone{
In fact, a P(EC)\textsuperscript{$M$}, $M > 1$,
implementation would collect and put aside
the values $z_n^{(1)}, \dotsc, z_n^{(M-1)}$, which
seems unintuitive from an inference point of view,
where it is natural to assume that more data should
yield a better approximation. A natural question
would be whether this is a case of diminishing
returns of approximation quality for computational power,
but this is beyond the scope of this paper.
}

\begin{figure}
  \centering
  \setlength{\figurewidth}{0.45\textwidth}
  \setlength{\figureheight}{\figurewidth}
  \input{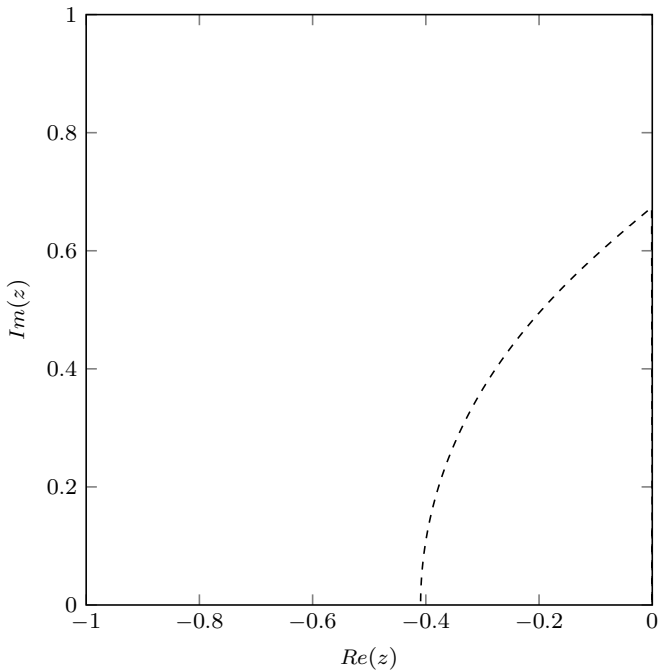}
  \caption{Partial \editone{stability domain} of the probabilistic Nordsieck method
    using the IWP($2$) in the negative real, positive imaginary quadrant.
    The method converges for step sizes $h$ on linear problems $y' = \lambda y$,
    if $h\lambda \ce z \in \mathbb{C}$ lies in the region of stability in the lower right corner.
    See \cite{DeuflhardScientific} for details.
  }
  \label{fig:root-locus-N2}
\end{figure}

One current limitation of the \editone{PFOS} is
its fixed integration order $q$ over the whole integration
domain $\mathbb{T}$. The reason for this limitation is that it is
conceptually not straight forward to connect spline models of different order at knots $\tn$.
However, the ability to adapt the integration order
during runtime has been key in improving the efficiency of
modern solvers \cite{byrne1975polyalgorithm}. Furthermore, the
method corresponding to the IWP($2$) model has a rather small
region of stability which is depicted in Figure~\ref{fig:root-locus-N2},
specially in comparison to \editone{backward differentiation formulas (BDFs)} \cite{DeuflhardScientific}.
This makes the method impractical for stiff equations.

\editone{
It is natural to ask what happens in the case
of the IWP($q$), $q > 2$. Using techniques from
the analysis of Kalman filters, one can show that
these models also contain a stable subsystem and
that the weights $\mat{K}_n$ will converge to
a fixed point $\mat{K}^*$, even for non-zero, but
constant, $R^2$. However,
}
it remains unclear whether
they will be practical. In particular, these methods might
even be unstable for most spline models \cite{Loscalzo1967Spline}.
We have tested the IWP($q$), $q \in \{1, \dotsc, 4\}$, empirically
on the Hull et al. benchmark (see Section~\ref{sec:experiments})
and have observed that these converge in practice on these
non-stiff problems.



\subsection{Initialization via Runge--Kutta methods}
\label{sec:connection-rk}

\editone{
Thus far, we have provided the analysis of the long-term
behavior of the algorithm, when several Kalman filter steps
have been computed and the steady-state is reached.
Crucially, a necessary condition for this analysis is
that enough updates have been performed such that the
observable space spans the entire state-space, which is
$q+1$ updates in the case for the IWP($q$).
}

\editone{
Thus, the question remains how to initialize the filter.
Schober et al.~\cite{schober2014nips} have shown that
there are Runge--Kutta steps that coincide with the \edit{maximum a posteriori (MAP)}
of the IWP($q$) for $q \leq 3$. This requires $q+1$
updates using a diffuse prior with
$\mat{C}_\tmo = \lim_{\mathcal{H} \to \infty} \mat{Q}(\mathcal{H})$.
In practice, one takes a Runge--Kutta step with the corresponding
formula and plugs the resulting values into the analytic
expressions for $\vec{m}_{t_1}$ and $\mat{C}_{t_1}$ at $t_1$.
Additionally to the cases presented in \cite{schober2014nips},
we can report a match between a four step Runge--Kutta formula
of order four and the IWP($4$). This match is obtained for
the evaluation knots $t_0 + c_ih$ with the vector
$\vec{c} = (0, \nicefrac{1}{3}, \nicefrac{1}{2}, 1)\Trans$.
Details and exact expressions are given in Appendix~\ref{sec:rk-derivations}.
This approach is structurally similar to an algorithm
given by Gear~\cite{gear1980rkstarters} for the case
of classical Runge--Kutta and Nordsieck methods.
}

\editone{
However, we want to stress that the analysis by Schober
et al.~\cite{schober2014nips} is done under 
exactly the same model and with the same assumptions that
have been presented here in different notation.
Therefore, the initialization does not require a separate
model and our requirement of a globally-defined solver
still holds.
}

\editone{
Finally, it should be pointed out that this is only one
feasible initialization.
In cases where automatic differentiation \cite{Griewank2008EDP}
is available, this can be used to initialize the Nordsieck
vector up to numerical precision and set $\mat{C}_\tmo$ to $\mat{0}$.
Nordsieck originally proposed \cite{nordsieck1962numerical}
start with an initial vector
$\vec{m}_\tmo = 0$, followed by $q+1$ steps with positive and
$q+1$ with negative direction (that is, integrating backwards to the start).
One interpretation is that the method uses
$\vec{m}_\tmo = \Exp[X_\tmo\g \tilde{z}_{-1}, \dots, \tilde{z}_q]$,
with tentative $\tilde{z}_i$ computed out of this process.
}




\section{Error estimation and hyperparameter adaptation}
\label{sec:error-estimation}

While the general algorithm described in Sec.~\ref{sec:nordsieck}
can be applied to any IVP at this stage, a modern ODE solver also
requires the ability to automatically select sensible values for its
hyper-parameters.
The filter has three remaining parameters to choose:
the dimensionality $q$ of the state space, the diffusion amplitude $\ssq$
and the step size $h$.

To obtain a globally consistent probability distribution, we fix $q = 2$
throughout the integration to test the third order method presented in Sec.~\ref{sec:nordsieck}.
For the remaining two parameters, we first note that estimating $\ssq$
will lend itself naturally to choose the step size.
To see this, one needs to make the connection to classical ODE solvers and the
interpretation of the state space model. In classical ODE solvers, $h_n$
is determined based on local error analysis, that is, $h_n$ is a function
of the error $e_\tn$ introduced from step $\tnmo$ to step $\tn$.
Then, $h_n$ is computed as a function of the allowed tolerance and the
expected error which is assumed to evolve similarly to the current error.

As is common in solving IVPs, we base error estimation on \emph{local errors}.
Assume that the predicted solution $\vec{m}_\tnmo$ at time $\tnmo$ is
error free, that is, $\mat{C}_\tnmo = \mat{0}$. Then, by
Equations~\eqref{eq:pred-cov} and \eqref{eq:kalman-res}, we have
\begin{equation}
p(\lambda_n\g\ssq) 
= \N (\lambda_n;  z_n - \mat{H}_1\vec{m}_\tn^-,
      \mat{H}_1\ssq\mat{\bar{Q}}(h)\mat{H}_1\Trans).
\label{eq:res-likelihood}
\end{equation}
One way to find the optimal $\ssq$ is to construct the maximum likelihood
estimator from Equation~\eqref{eq:res-likelihood} which is given by
\begin{equation}
\begin{aligned}
\hat{\sigma}^2 &= 
(z_n - \mat{H}_1\vec{m}_\tnmo^-)\Trans
(\mat{H}_1\mat{\bar{Q}}(h)\mat{H}_1\Trans)^{-1}
(z_n - \mat{H}_1\vec{m}_\tnmo^-)\\
  &= \frac{(z_n - \mat{H}_1\vec{m}_\tnmo^-)^2}{\mat{H}_1\mat{\bar{Q}}(h)\mat{H}_1\Trans}.
\end{aligned}
\label{eq:ml-estimator}
\end{equation}
For the last equation we used the fact that all the involved quantities are scalars.

To allow for a greater flexibility of the model, we allow amplitude $\ssq$ to
vary for different steps $\ssq_{t_n}$. Note, that the mean values are then no
longer independent of $\ssq$, because the factor no longer cancels out in the
computation of $K_n$ in Equation~\eqref{eq:update-mean}. However, this situation
is indeed intended: if there was more diffusion in $[t_{n-1}, t_n]$, we want a
stronger update to the mean solution as the observed value is more informative.
\edit{Additionally, Equation~\eqref{eq:kalman-res}}
is independent of $\ssq_{t_n}$ or any other covariance information $\mat{P}_{t_n}^-, \mat{Q}(h)$.
Therefore, we can apply Equation~\eqref{eq:kalman-res} before Equation~\eqref{eq:pred-cov},
update $\ssq_{t_n}$ and then continue to compute the rest of the Kalman step.
This idea is similar in spirit to \cite[\S 11]{jazwinski1970stochastic},
but follows the general idea of error estimation in numerical ODE solvers,
where local error information is available only.

At this point, the
inference interpretation of numerical computation comes to bear:
once the initial modeling decision---modeling a deterministic object with
a probability measure to describe the uncertainty over the solution---is accepted,
everything else follows naturally from the probabilistic description. Most
importantly, there are no neglected higher-order terms, as they are all
incorporated in the diffusion assumption.

This kind of lightweight error estimation is 
a key ingredient to probabilistic numerical methods: one goal of a
probabilistic model are improved decisions under \emph{uncertainty}.
This uncertainty is necessarily a crude approximation, since a more
accurate error estimator could be used to improve the overall solution
quality. However, the reduction of computational efforts up to a
tolerated error is exactly what modern numerical solvers try to achieve.


This error estimate can now be used in the conventional way of adapting the
step size which we will restate here to give a complete description of the
inference algorithm (see also~\cite{byrne1975polyalgorithm}).
Given an error weighting vector $\vec{w}$,
the algorithm computes the weighted expected error
\begin{equation}
(\vec{D}_\tn)_i = (\mat{H_1} \ssq_\tn \mat{\bar{Q}}(h_n) \mat{H_1}\Trans)_i^{\edit{1/2}} \vec{w}_i,
\end{equation}
where $\mat{\bar{Q}}(h_n) = [\ssq_\tn]^{-1} \mat{Q}(h_n)$ is the normalized
diffusion matrix,
and checks whether some error tolerance with parameter $\epsilon$ is met
\newcommand{\tol}{\epsilon}
\begin{equation}
\vec{D}_\tn \leq \bar{\tol} \ce \tol \frac{h_n}{S} \label{eq:error-test}
\end{equation}
where $h_n$ is the step length and $S$ can be either chosen to be $S = 1$
(error per unit step) or $S = h_n$ (error per step).
If Equation~\eqref{eq:error-test} holds, the step is accepted and integration
continues. Otherwise, the step is rejected as too inaccurate and the step
is repeated. In both cases, a new step length is computed which will likely
satisfy Equation~\eqref{eq:error-test} on the next step attempt. The new step size
is computed as
\begin{equation}
h_{n+1} = \rho \left(\frac{\bar{\tol}}{\vec{D}_{t_n}}\right)^{\frac{1}{q+1}}
\end{equation}
where $\rho \in (0, 1),\,\rho \approx 1$ is a safety factor. Additionally,
we also follow \del{established} best practices \cite{hairer87:_solvin_ordin_differ_equat_i} limit
the rate of change $\eta_{\min} < \nicefrac{h_{n+1}}{h_n} < \eta_{\max}$.
In our code, we set $\rho \ce 0.95,\;\eta_{\min} \ce 0.1$ and $\eta_{\max} \ce 5$.

\subsection{Global versus local error estimation}
\label{sec:global-local}

\edit{
The results presented in preceding sections pertain to the estimation
of \emph{local} extrapolation errors.
It is a well known aspect of ODE solvers
\cite[\textsection III.5]{hairer87:_solvin_ordin_differ_equat_i}
that the \emph{global} error can be exponentially larger
than the local error.
More precisely, to scale the stochastic process such that the variance
of the resulting posterior measure relates to the square \emph{global} error,
the intensity $\ssq_n$ of the stochastic process must be multiplied by a factor
\cite[Thm III.5.8]{hairer87:_solvin_ordin_differ_equat_i}
$\exp(L^*(T-t_0))$, where $L^*$ is a constant depending on
the problem.
Although related, $L^*$ is not the same as the local Lipschitz
constant $L$, and harder to estimate in practice (more details in \cite[\textsection III.5]{hairer87:_solvin_ordin_differ_equat_i}). We stress that this issue does not invalidate the probabilistic interpretation of the posterior measure as such. It is just that the scale of the posterior has to be estimated differently if the posterior is supposed to capture global error instead of local error. In practice, the global error estimate resulting from this re-scaling is often very conservative.
}


\section{Experiments}
\label{sec:experiments}

\edit{
To evaluate the model, we provide
two sets of experiments.
First, we qualitatively examine the
uncertainty quantification by visualizing the
posterior distribution of two example problems.
We also compare our proposed observation assumption
against the model described by Chkrebtii et al.~\cite{o.13:_bayes_uncer_quant_differ_equat}.
Second, we more rigorously evaluate the solver on a
benchmark, and compare it to existing non-probabilistic codes. Our goal in this work is to construct an algorithm that produces meaningful probabilistic output at a computational cost that is comparable to that of classic, non-probabilistic solvers. The experiments will show that this is indeed possible. Other probabilistic methods, in particular that of Chkrebtii et al.~\cite{o.13:_bayes_uncer_quant_differ_equat}, aim for a more expressive, non-Gaussian posterior. In exchange, the computational cost of these methods is at least a large multiple of that of the method proposed here, or even polynomially larger. These methods and ours differ in their intended use-cases: More elaborate but expensive posteriors are valuable for tasks in which uncertainty quantification is a central goal, while our solver aims to provide a meaningful posterior as additional functionality in settings where fast estimates are the principal objective.
}

\subsection{Uncertainty quantification}
\label{sec:visualization}

\edit{
We apply the probabilistic filtering ODE solver
on two problems with attracting orbits:
the Brusselator \cite{Levefer1971Chemical}
and van der Pol's equation \cite{vdP1926relaxation}.
The filter is applied twice on each problem,
once with a fixed step size and once with the
adaptive step size algorithm described in
Sec.~\ref{sec:error-estimation}.
To get a visually interesting plot, the fixed step size and the tolerance threshold
were chosen as large as possible without causing instability.
Both cases are modeled with a local diffusion
parameter $\ssq_n$ which is estimated using
the maximum likelihood estimator of Sec.~\ref{sec:error-estimation}.
In the following plots, the samples use the
scale $\ssq_n$ arising from the local error estimate.
Because these systems are attractive, the global
error correction mentioned in Sec.~\ref{sec:global-local}
would lead to significantly more conservative uncertainty.
}

\edit{
The Brusselator is the idealized and simplified model
of an autocatalytic multi-molecular chemical reaction \cite{Levefer1971Chemical}.
The rate equations for the oscillating reactants are
\begin{equation}
	\begin{aligned}
		y_1' &= A + y_1^2y_2 - (B+1)y_1\\
		y_2' &= By_1 - y_1^2y_2,
	\end{aligned}
	\label{eq:brusselator}
\end{equation}
where $A$ and $B$ are positive constants describing
the initial concentrations of two reactants.
Following \cite{hairer87:_solvin_ordin_differ_equat_i},
we set $A = 1$, $B = 3$ and $(y_1(0), y_2(0))\Trans = (1.5, 3)\Trans$.
The integration domain $\mathbb{T} = [0, 10]$ has
been chosen such that the solution completes one
cycle on the attractor after an initial convergence
phase.
}

\begin{figure*}
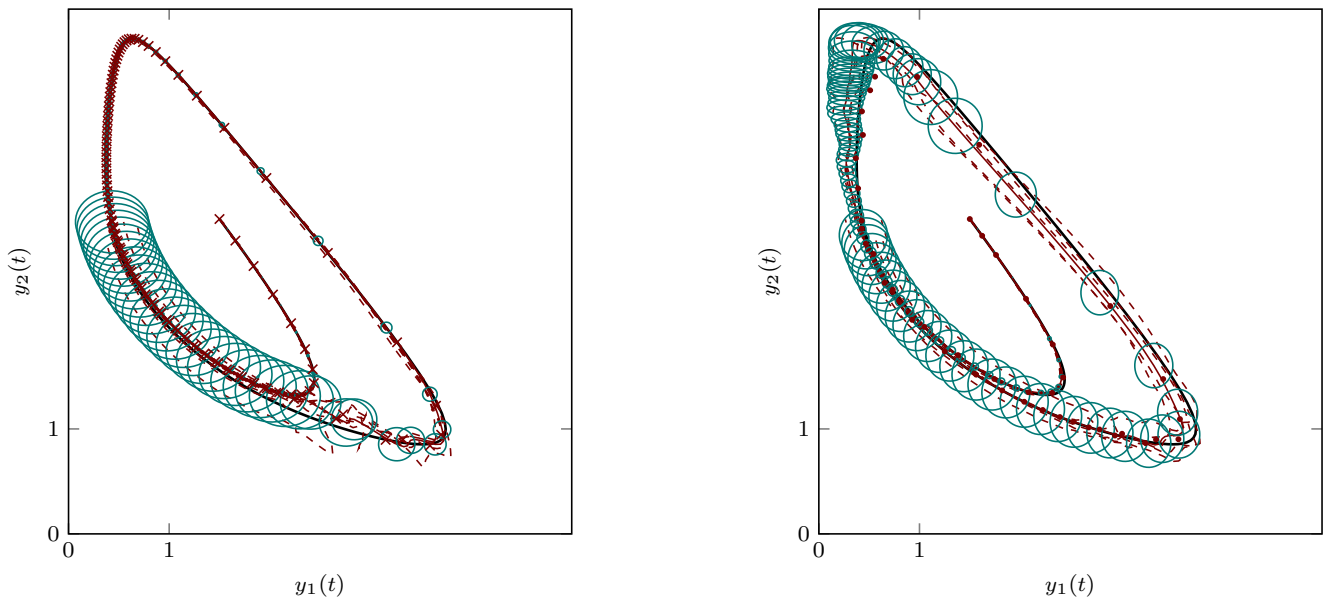

  \centering
  \setlength{\figurewidth}{0.4\textwidth}
  \setlength{\figureheight}{\figurewidth}
  \input{BRUSWOD-fixed}
  \hfill
  \input{BRUSWOD-adapt}
  \caption{Numerical solution of the Brusselator~\eqref{eq:brusselator}
    using the probabilistic filtering ODE solver.
    The plots show the solution computed by \texttt{ode45} using
    $\texttt{RelTol} = \texttt{AbsTol} = 1 \times 10^{-13}$
    (black, background), the posterior mean (red, thick line),
    iso-contourlines of twice the posterior standard deviation at a subsample
    of the knots (green) and samples from the posterior distribution (red, dashed lines).
    \textbf{Left:} Using a fixed step size of $h = h_n = 0.0834$.
    The computation requires $120$ steps.
    \textbf{Right:} Using the adaptive step size selection with
    error weighting $w_i(y) = (\tau y_i + \tau)^{-1}, \tau = 0.1$.
    The computation requires $43$ steps.
    See \cite[\textsection 1.6]{hairer87:_solvin_ordin_differ_equat_i} for details.
  }
  \label{fig:bruswod}
\end{figure*}

\edit{
The results in Figure~\ref{fig:bruswod} demonstrate
the effectiveness of the error estimator.
This problem also demonstrates the quality and
utility of the step size adaptation algorithm,
since on the majority of the solution trajectory
the algorithm is not limited by stability constrains.
}

\edit{
Van der Pol's equation \cite{vdP1926relaxation}
describes an oscillation with a non-linear damping
factor $\alpha$
\begin{equation}
	\begin{aligned}
		0 &= y'' + \alpha y' + y\\
		\alpha &= \mu (y^2 - 1)
	\end{aligned}
	\label{eq:vdp}
\end{equation}
with a positive constant $\mu > 0$.
Originally, this model has been used to describe
vacuum tube curcuits.
The limit cycle alternates between a non-stiff
phase of rapid change and a stiff phase of slow decay.
The larger $\mu$ the more pronounced both effects are.
In our example, we set $\mu = 1$ and integrate over
one period with the initial value on the graph of
the limit cycle.
Exact values can be found in
\cite[\textsection I.16]{hairer87:_solvin_ordin_differ_equat_i}.
}

\begin{figure*}
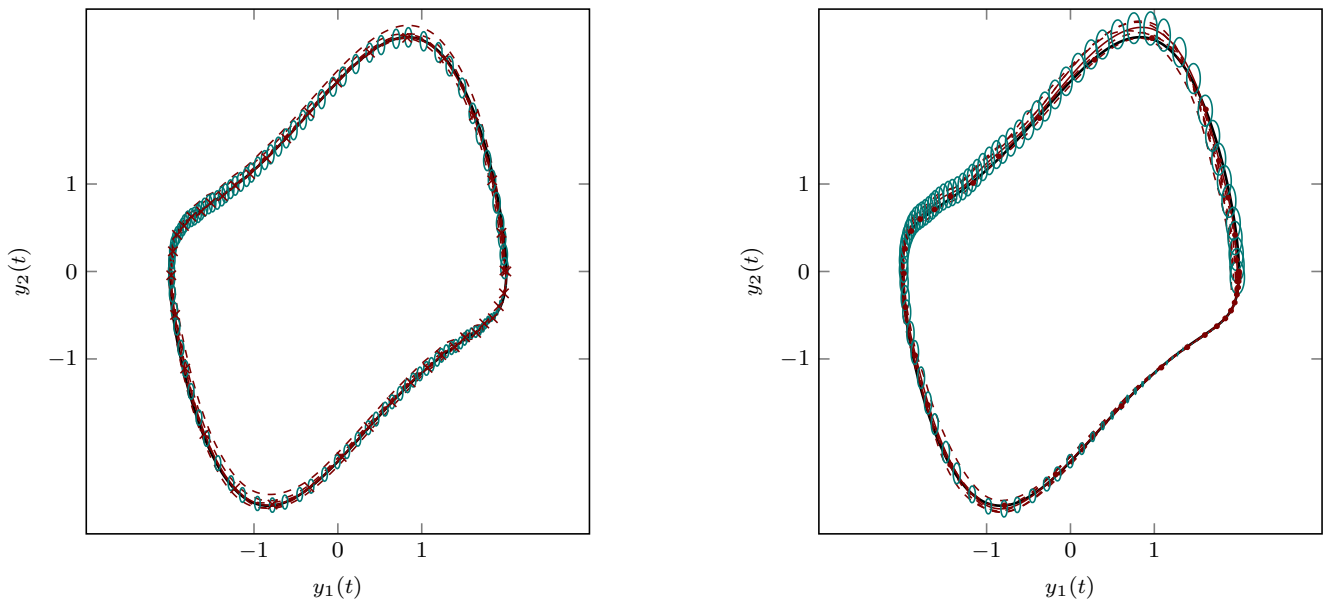

  \centering
  \setlength{\figurewidth}{0.4\textwidth}
  \setlength{\figureheight}{\figurewidth}
  \input{VDPOL-fixed}
  \hfill
  \input{VDPOL-adapt}
  \caption{
    Numerical solution of van der Pol's equation~\eqref{eq:vdp}
    using the probabilistic filtering ODE solver, integrated over
    one limit cycle period $\mathbb{T} = [0, T]$ with initial
    value $y(0) = (A, 0)\Trans$, where $T \approx 6.6633$ and $A \approx 2.0086$.
    The plots show the solution computed by \texttt{ode45} using
    $\texttt{RelTol} = \texttt{AbsTol} = 1 \times 10^{-13}$
    (black, background), the posterior mean (red, thick line),
    iso-contourlines of twice the posterior standard deviation at a subsample
    of the knots (green) and samples from the posterior distribution (red, dashed lines).
    \textbf{Left:} Using a fixed step size of $h = h_n = 0.1667$.
    The computation requires $40$ steps.
    \textbf{Right:} Using the adaptive step size selection with
    error weighting $w_i(y) = (\tau y_i + \tau)^{-1}, \tau = 0.1$.
    The computation requires $41$ steps.
    See \cite[\textsection 1.6]{hairer87:_solvin_ordin_differ_equat_i} for details.
  }
  \label{fig:vdpol}
\end{figure*}

\edit{
Figure~\ref{fig:vdpol} plots the filter results.
In the case of van der Pol's equation, the benefit
of step size adaptation is essentially nil, because conservative adaptation consumes the gains
on the non-stiff parts. However, the example demonstrates the capability
to learn an anisotropic diffusion model for individual components.}

\begin{figure*}
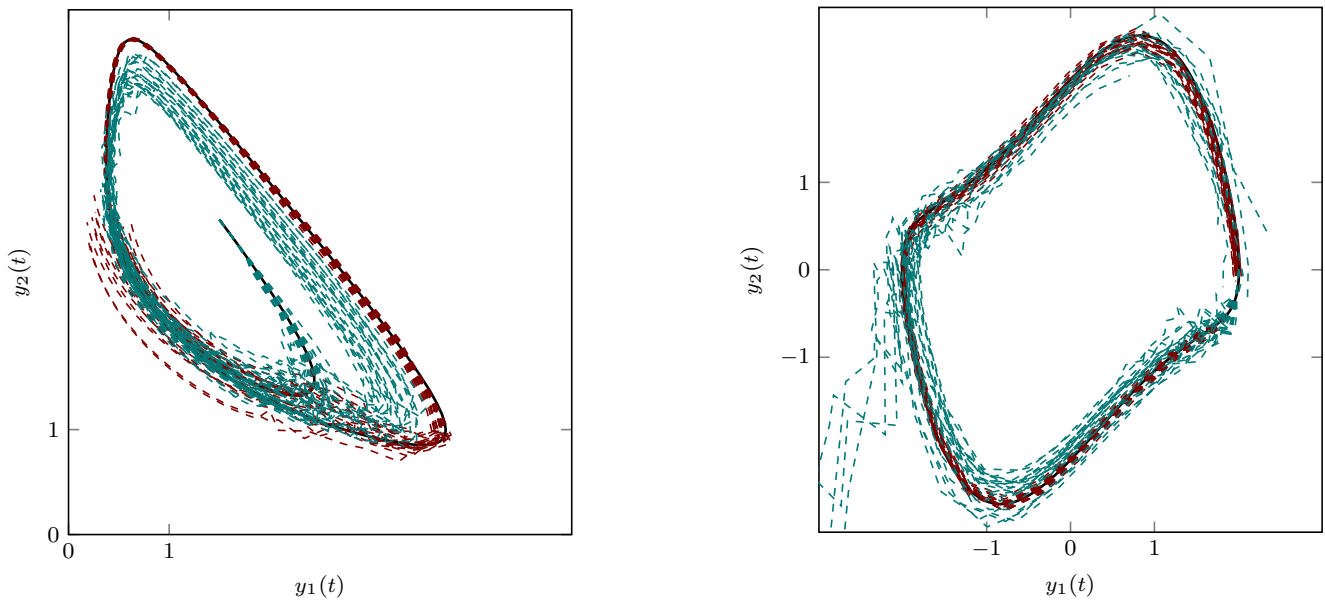

  \centering
  \setlength{\figurewidth}{0.4\textwidth}
  \setlength{\figureheight}{\figurewidth}
  \input{BRUSWOD-sampled}
  \hfill
  \input{VDPOL-sampled}
  \caption{
  	Comparison of two different evaluation strategies on problems~\eqref{eq:brusselator}
  	and \eqref{eq:vdp}.
  	\textcolor{dred}{Red:} samples from the posterior
  	as in Figs.~\ref{fig:bruswod} and \ref{fig:vdpol}.
  	\textcolor{MPG}{Green:} Similar, but evaluating at
  	$z_n = f(t_n, (u_\tn)_0), u_\tn \sim \N(\vec{m}^-_\tn, \mat{C}^-_\tn)$.
  	This is similar to \cite{o.13:_bayes_uncer_quant_differ_equat}.
  }
  \label{fig:comp-chkr}
\end{figure*}

\edit{
Finally, we compare two different strategies of quantifying
the uncertainty.
To this end, we compare our proposed model to the observation
model proposed by Chkrebtii et al.~\cite[\textsection 3.1]{o.13:_bayes_uncer_quant_differ_equat}.
In this case, we set $z_n = f(\tn, (u_\tn)_0),  u_\tn \sim \N(\vec{m}^-_\tn, \mat{C}^-_\tn)$.
Figure~\ref{fig:comp-chkr} shows samples of
the posterior distribution,
computed with two different evaluation schemes. This scheme is not exactly the same as the one proposed by Chkrebtii et al---their algorithm actually has cubic complexity in the number of $f$-evaluations, thus is  limited to a relatively small number of evaluation steps. But our version captures the principal difference between their algorithm and the simpler filter proposed here: Their algorithm draws separate samples involving independent evaluations of $f$ at perturbed locations, while ours draws samples from a single posterior constructed from one single set of $f$-evaluations. 
As expected, the model of Chkrebtii et al.\ provides
a richer output structure, for example, by identifying
divergent solutions (right subplot) if the solver
leaves the region of attraction.
However, to obtain individual samples, the entire
algorithm has to run repeatedly, so the cost of producing $S$ samples is $S$ times that of our algorithm, which produces all its samples in one run, without requiring additional evaluations of $f$.
}



\subsection{Benchmark evaluation}
\label{sec:benchmark}

\newcommand{\detest}{\texttt{DETEST}}

As is the case with many modern solvers,
the theoretical guarantees do not extend
to the full implementation with error
estimation and step size control.
Therefore, an empirical assessment is
necessary to compare against trusted
implementations.
We compare the proposed Kalman filter to
a representative set of standard algorithms
on the \detest{} benchmark set \cite{hull1972comparing}.
While other standardized tests have been proposed
\cite{crane1969comparative,krogh1973testing},
\detest{} has repeatedly been described as
representative \cite{shampine1976solving,deuflhard1983order}.
By choosing the same comparison criteria across
all test problems and tested implementations,
the benchmark provides the necessary data to
make predictions on the behavior on a large
class of problems.

Two different dimensions of performance are
considered in \cite{hull1972comparing}: the
\emph{computational cost} and the \emph{solution
quality}. Computational cost is reported in
execution time (in seconds) and number of
function evaluations (abbreviated as \#FE). Although the former
is more relevant in practice, we only report
the latter here as the codes in
\cite{hull1972comparing} are outdated and our
proof-of-concept code is not yet optimized for speed.
Nevertheless, since the execution times are
proportional to the \#FE, this provides a
reliable estimator of computational efficiency.
\detest{} only considers methods with automatic
step size adaptation, and thus measures the
solution quality by comparing the local error
with the requested tolerance $\tol$. A code
is considered to produce high quality solutions
if the results are within the requested tolerance,
but are also not of excessive unrequested higher
accuracy. Therefore, errors are reported per
unit step.
Reported are the maximum error
$\max \{\xi_n[h_n \tol]^{-1}\g n = 1, \dotsc, N\}$
per unit step and the
percentage of deceived steps
$|\{\xi_n\g \xi_n > h_n \tol, n = 1,\dotsc, N\}|/N$,
where the local errors $\xi_n$ are defined as
$\norm{y_\tn - y(\tn; y(\tnmo) = y_\tnmo)}_\infty$ and
$y(\tn; y(\tnmo) = y_\tnmo)$ defines the IVP
$y' = f(t, y),\,y(\tnmo) = y_\tnmo,\, t \in [\tnmo, \tn]$.

\begin{table*}[tbp]
\caption{Summary of \texttt{DETEST} results}
\label{tab:detest-summary}

\begin{center}\footnotesize
\begin{tabular}{l*{3}{r}}
\toprule
Method & Total fcn. evals. & Avg. \% deceived & Max. error \\
\midrule
\noalign{\smallskip}
$\tol = 10^{-3}$ \\
\cmidrule{1-1}
\noalign{\smallskip}
Extrapolation      & 16553  &  2.0 &    7.8 \\
Adams (Krogh)      &  5394  &  1.1 &    5.3 \\
Adams (Gear)       &  9498  &  0.9 &    1.5 \\
RK (4th, Kutta)    &  8363  &  5.1 &   25.9 \\
RK (6th, Butcher)  & 11105  &  5.1 & 1788.1 \\
RK (8th, Shanks)   & 12355  &  6.3 & 1120.6 \\
RK (3th, Shampine) & 15085  &  5.9 &    2.4 \\
RK (5th, Shampine) &  5785  & 11.2 &    9.5 \\
Adams (Shampine)   &  5692  &  6.5 &    7.7 \\
PNM                & 19091  &  0.2 &    1.5 \\
\midrule
\noalign{\smallskip}
$\tol = 10^{-6}$ \\
\cmidrule{1-1}
\noalign{\smallskip}
Extrapolation      &  26704 &  0.1 &   2.3 \\
Adams (Krogh)      &  11353 &  1.4 &   7.3 \\
Adams (Gear)       &  18155 &  0.8 &   2.6 \\
RK (4th, Kutta)    &  30763 &  1.8 &  29.1 \\
RK (6th, Butcher)  &  23540 &  1.6 & 142.5 \\
RK (8th, Shanks)   &  20493 &  4.2 &   4.7 \\
RK (3th, Shampine) & 430975 &  0.0 &   1.9 \\
RK (5th, Shampine) &  19879 &  0.0 &   1.1 \\
Adams (Shampine)   &  10777 &  3.6 &   6.3 \\
PNM                & 405469 &  0.0 &   1.4 \\
\midrule
\noalign{\smallskip}
$\tol = 10^{-9}$ \\
\cmidrule{1-1}
\noalign{\smallskip}
Extrapolation      &    43054 & 0.0 &    0.6 \\
Adams (Krogh)      &    18984 & 0.5 &    4.0 \\
Adams (Gear)       &    38439 & 2.3 &    2.7 \\
RK (4th, Kutta)    &   146262 & 0.3 &    2.9 \\
RK (6th, Butcher)  &    58634 & 0.9 &  443.4 \\
RK (8th, Shanks)   &    39663 & 2.1 &   20.9 \\
RK (3th, Shampine) & 13587187 & 3.1 &  689.0 \\
RK (5th, Shampine) &   103345 & 0.1 &    2.4 \\
Adams (Shampine)   &    18274 & 2.2 &   11.5 \\
PNM                & 12731730 & 4.5 & 1938.0 \\
\bottomrule
\end{tabular}
\end{center}

\end{table*}

Here, we report the results from the proposed solver
originating from the IWP($2$) model
as well as the results from the original
Hull et al.\ paper \cite{hull1972comparing}.
We have not been able to obtain a copy of the
codes used in Hull et al.\ and only report
their numbers for sake of completeness.
We also ran the tests on the solvers provided in \matlab{}.
Table~\ref{tab:detest-summary} lists the summary results
for all methods 
and all tolerances.
For detailed results on individual problems see
Figures~\ref{fig:fcnCalls}--\ref{fig:maxError} in Sec.~\ref{sec:detailed-results}.
For a complete and detailed description of the
benchmark, we refer to \cite{hull1972comparing}.
Our implementation is publicly available.\repo{}

\begin{figure*}
  \centering
  \setlength{\figurewidth}{0.165\textwidth}
  \setlength{\figureheight}{\figurewidth}
%
%
\definecolor{mycolor1}{rgb}{0.50000,0.50000,0.75450}%
\definecolor{mycolor2}{rgb}{0.25000,0.25000,0.63175}%
\definecolor{mycolor3}{rgb}{0.14645,0.14645,0.58091}%
\definecolor{mycolor4}{rgb}{0.06699,0.06699,0.54189}%
\definecolor{mycolor5}{rgb}{0.00000,0.00000,0.50900}%
\definecolor{mycolor6}{rgb}{0.00000,0.47170,0.46040}%
\definecolor{mycolor7}{rgb}{0.49060,0.00000,0.00000}%
\begin{tikzpicture}

\begin{axis}[%
width=0.951\figurewidth,
height=\figureheight,
at={(0\figurewidth,0\figureheight)},
scale only axis,
xmin=0,
xmax=1.6,
xtick={0,0.25,0.5,0.75,1,1.25,1.5,1.55},
xticklabels={{0},{},{0.5},{},{1},{},{},{$>1.5$}},
ymin=0,
ymax=1,
ytick={0,0.1,0.2,0.3,0.4,0.5,0.6,0.7,0.8,0.9,1},
yticklabels={\empty},
axis background/.style={fill=white},
mystyle
]
\addplot [color=mycolor1,solid,forget plot]
  table[row sep=crcr]{%
0	0.256637168141593\\
0.05	0.513274336283186\\
0.15	0.68141592920354\\
0.25	0.902654867256637\\
0.35	0.920353982300885\\
0.45	0.955752212389381\\
0.55	0.955752212389381\\
0.65	0.964601769911504\\
0.75	0.964601769911504\\
0.85	0.964601769911504\\
0.95	0.964601769911504\\
1.05	0.964601769911504\\
1.15	0.964601769911504\\
1.25	0.964601769911504\\
1.35	0.964601769911504\\
1.45	0.964601769911504\\
1.55	1\\
};
\addplot [color=mycolor2,solid,forget plot]
  table[row sep=crcr]{%
0	0.378378378378378\\
0.05	0.756756756756757\\
0.15	0.905405405405405\\
0.25	0.932432432432432\\
0.35	0.945945945945946\\
0.45	0.959459459459459\\
0.55	0.959459459459459\\
0.65	0.959459459459459\\
0.75	0.959459459459459\\
0.85	0.959459459459459\\
0.95	0.959459459459459\\
1.05	0.959459459459459\\
1.15	0.959459459459459\\
1.25	0.959459459459459\\
1.35	0.959459459459459\\
1.45	0.959459459459459\\
1.55	1\\
};
\addplot [color=mycolor3,solid,forget plot]
  table[row sep=crcr]{%
0	0.462986198243413\\
0.05	0.925972396486826\\
0.15	0.968632371392723\\
0.25	0.987452948557089\\
0.35	0.991217063989962\\
0.45	0.994981179422836\\
0.55	0.994981179422836\\
0.65	0.996235884567127\\
0.75	0.996235884567127\\
0.85	0.996235884567127\\
0.95	0.996235884567127\\
1.05	0.996235884567127\\
1.15	0.996235884567127\\
1.25	0.996235884567127\\
1.35	0.996235884567127\\
1.45	0.996235884567127\\
1.55	1\\
};
\addplot [color=mycolor4,solid,forget plot]
  table[row sep=crcr]{%
0	0.384615384615385\\
0.05	0.769230769230769\\
0.15	0.901098901098901\\
0.25	0.956043956043956\\
0.35	0.956043956043956\\
0.45	0.967032967032967\\
0.55	0.967032967032967\\
0.65	0.967032967032967\\
0.75	0.967032967032967\\
0.85	0.967032967032967\\
0.95	0.967032967032967\\
1.05	0.967032967032967\\
1.15	0.967032967032967\\
1.25	0.967032967032967\\
1.35	0.967032967032967\\
1.45	0.967032967032967\\
1.55	1\\
};
\addplot [color=mycolor5,solid,forget plot]
  table[row sep=crcr]{%
0	0.433333333333333\\
0.05	0.866666666666667\\
0.15	0.944444444444444\\
0.25	0.955555555555556\\
0.35	0.966666666666667\\
0.45	0.966666666666667\\
0.55	0.966666666666667\\
0.65	0.966666666666667\\
0.75	0.966666666666667\\
0.85	0.966666666666667\\
0.95	0.966666666666667\\
1.05	0.966666666666667\\
1.15	0.966666666666667\\
1.25	0.966666666666667\\
1.35	0.966666666666667\\
1.45	0.966666666666667\\
1.55	1\\
};
\addplot [color=mycolor6,dashed,forget plot]
  table[row sep=crcr]{%
0	0\\
0.0161616161616162	0.0128945426727563\\
0.0323232323232323	0.0257857178319707\\
0.0484848484848485	0.0386701606022311\\
0.0646464646464646	0.0515445113809416\\
0.0808080808080808	0.0644054184661257\\
0.096969696969697	0.0772495406739226\\
0.113131313131313	0.0900735499423687\\
0.129292929292929	0.102874133918085\\
0.145454545454545	0.115647998522515\\
0.161616161616162	0.128391870494411\\
0.177777777777778	0.141102499905278\\
0.193939393939394	0.153776662644587\\
0.21010101010101	0.166411162871553\\
0.226262626262626	0.179002835430406\\
0.242424242424242	0.191548548226099\\
0.258585858585859	0.204045204557483\\
0.274747474747475	0.216489745405068\\
0.290909090909091	0.228879151670556\\
0.307070707070707	0.241210446365418\\
0.323232323232323	0.253480696745887\\
0.339393939393939	0.265687016391836\\
0.355555555555556	0.277826567227089\\
0.371717171717172	0.289896561478844\\
0.387878787878788	0.30189426357398\\
0.404040404040404	0.313816991970125\\
0.42020202020202	0.325662120919495\\
0.436363636363636	0.337427082163601\\
0.452525252525253	0.349109366557089\\
0.468686868686869	0.360706525619037\\
0.484848484848485	0.372216173010211\\
0.501010101010101	0.383635985934891\\
0.517171717171717	0.394963706465981\\
0.533333333333333	0.406197142792298\\
0.54949494949495	0.417334170387006\\
0.565656565656566	0.428372733096346\\
0.581818181818182	0.439310844147906\\
0.597979797979798	0.450146587077822\\
0.614141414141414	0.460878116576438\\
0.63030303030303	0.47150365925208\\
0.646464646464647	0.482021514312714\\
0.662626262626263	0.492430054165421\\
0.678787878787879	0.502727724933718\\
0.694949494949495	0.512913046892906\\
0.711111111111111	0.522984614823739\\
0.727272727272727	0.532941098284826\\
0.743434343434344	0.542781241804333\\
0.75959595959596	0.552503864991613\\
0.775757575757576	0.562107862569582\\
0.791919191919192	0.571592204328704\\
0.808080808080808	0.580955935003612\\
0.824242424242424	0.590198174073463\\
0.84040404040404	0.599318115487263\\
0.856565656565657	0.608315027315463\\
0.872727272727273	0.617188251329263\\
0.888888888888889	0.625937202509117\\
0.905050505050505	0.634561368484056\\
0.921212121212121	0.643060308903506\\
0.937373737373738	0.651433654743375\\
0.953535353535354	0.65968110754825\\
0.96969696969697	0.667802438611626\\
0.985858585858586	0.675797488096148\\
1.0020202020202	0.683666164095923\\
1.01818181818182	0.691408441643\\
1.03434343434343	0.699024361660193\\
1.05050505050505	0.706514029862441\\
1.06666666666667	0.713877615608982\\
1.08282828282828	0.721115350708616\\
1.0989898989899	0.728227528180404\\
1.11515151515152	0.73521450097215\\
1.13131313131313	0.742076680639059\\
1.14747474747475	0.748814535984982\\
1.16363636363636	0.75542859166865\\
1.17979797979798	0.761919426777358\\
1.1959595959596	0.768287673370511\\
1.21212121212121	0.774534014995493\\
1.22828282828283	0.780659185178288\\
1.24444444444444	0.786663965891287\\
1.26060606060606	0.792549186000707\\
1.27676767676768	0.79831571969601\\
1.29292929292929	0.803964484903738\\
1.30909090909091	0.809496441688096\\
1.32525252525253	0.814912590640635\\
1.34141414141414	0.820213971261342\\
1.35757575757576	0.8254016603334\\
1.37373737373737	0.830476770293846\\
1.38989898989899	0.835440447602338\\
1.40606060606061	0.840293871110162\\
1.42222222222222	0.845038250431579\\
1.43838383838384	0.849674824319579\\
1.45454545454545	0.854204859048024\\
1.47070707070707	0.858629646802131\\
1.48686868686869	0.862950504079186\\
1.5030303030303	0.867168770101306\\
1.51919191919192	0.871285805242022\\
1.53535353535354	0.875302989468395\\
1.55151515151515	0.879221720800282\\
1.56767676767677	0.883043413788345\\
1.58383838383838	0.886769498012308\\
1.6	0.890401416600884\\
};
\addplot [color=mycolor7,solid,forget plot]
  table[row sep=crcr]{%
1	0\\
1	1\\
};
\end{axis}
\end{tikzpicture}%
%
%
\definecolor{mycolor1}{rgb}{0.50000,0.50000,0.75450}%
\definecolor{mycolor2}{rgb}{0.25000,0.25000,0.63175}%
\definecolor{mycolor3}{rgb}{0.14645,0.14645,0.58091}%
\definecolor{mycolor4}{rgb}{0.06699,0.06699,0.54189}%
\definecolor{mycolor5}{rgb}{0.00000,0.00000,0.50900}%
\definecolor{mycolor6}{rgb}{0.00000,0.47170,0.46040}%
\definecolor{mycolor7}{rgb}{0.49060,0.00000,0.00000}%
\begin{tikzpicture}

\begin{axis}[%
width=0.951\figurewidth,
height=\figureheight,
at={(0\figurewidth,0\figureheight)},
scale only axis,
xmin=0,
xmax=1.6,
xtick={0,0.25,0.5,0.75,1,1.25,1.5,1.55},
xticklabels={{0},{},{0.5},{},{1},{},{},{$>1.5$}},
ymin=0,
ymax=1,
ytick={0,0.1,0.2,0.3,0.4,0.5,0.6,0.7,0.8,0.9,1},
yticklabels={\empty},
axis background/.style={fill=white},
mystyle
]
\addplot [color=mycolor1,solid,forget plot]
  table[row sep=crcr]{%
0	0.464761443448777\\
0.05	0.929522886897554\\
0.15	0.978445144102688\\
0.25	0.988617098571083\\
0.35	0.991765560668443\\
0.45	0.993945265197384\\
0.55	0.995640590942117\\
0.65	0.996124969726326\\
0.75	0.997335916686849\\
0.85	0.997578106078954\\
0.95	0.997578106078954\\
1.05	0.997578106078954\\
1.15	0.997820295471058\\
1.25	0.998062484863163\\
1.35	0.998304674255268\\
1.45	0.998789053039477\\
1.55	1\\
};
\addplot [color=mycolor2,solid,forget plot]
  table[row sep=crcr]{%
0	0.22375\\
0.05	0.4475\\
0.15	0.565\\
0.25	0.915\\
0.35	0.96375\\
0.45	0.9775\\
0.55	0.9825\\
0.65	0.985\\
0.75	0.9875\\
0.85	0.98875\\
0.95	0.9925\\
1.05	0.9925\\
1.15	0.9925\\
1.25	0.9925\\
1.35	0.9925\\
1.45	0.9925\\
1.55	1\\
};
\addplot [color=mycolor3,solid,forget plot]
  table[row sep=crcr]{%
0	0.252032520325203\\
0.05	0.504065040650406\\
0.15	0.669376693766938\\
0.25	0.880758807588076\\
0.35	0.932249322493225\\
0.45	0.943089430894309\\
0.55	0.953929539295393\\
0.65	0.956639566395664\\
0.75	0.959349593495935\\
0.85	0.959349593495935\\
0.95	0.962059620596206\\
1.05	0.967479674796748\\
1.15	0.967479674796748\\
1.25	0.967479674796748\\
1.35	0.967479674796748\\
1.45	0.967479674796748\\
1.55	1\\
};
\addplot [color=mycolor4,solid,forget plot]
  table[row sep=crcr]{%
0	0.446793002915452\\
0.05	0.893586005830904\\
0.15	0.951166180758018\\
0.25	0.967444120505345\\
0.35	0.977405247813411\\
0.45	0.982264334305151\\
0.55	0.987609329446064\\
0.65	0.990038872691934\\
0.75	0.992468415937804\\
0.85	0.993683187560739\\
0.95	0.9944120505345\\
1.05	0.996598639455782\\
1.15	0.996841593780369\\
1.25	0.997084548104956\\
1.35	0.99757045675413\\
1.45	0.998056365403304\\
1.55	1\\
};
\addplot [color=mycolor5,solid,forget plot]
  table[row sep=crcr]{%
0	0.391518131530424\\
0.05	0.783036263060848\\
0.15	0.888137676705593\\
0.25	0.926859250153657\\
0.35	0.952059004302397\\
0.45	0.964966195451752\\
0.55	0.972956361401352\\
0.65	0.979102642901045\\
0.75	0.982175783650891\\
0.85	0.985863552550707\\
0.95	0.987707437000615\\
1.05	0.988322065150584\\
1.15	0.99139520590043\\
1.25	0.992624462200369\\
1.35	0.993239090350338\\
1.45	0.993853718500307\\
1.55	1\\
};
\addplot [color=mycolor6,dashed,forget plot]
  table[row sep=crcr]{%
0	0\\
0.0161616161616162	0.0128945426727563\\
0.0323232323232323	0.0257857178319707\\
0.0484848484848485	0.0386701606022311\\
0.0646464646464646	0.0515445113809416\\
0.0808080808080808	0.0644054184661257\\
0.096969696969697	0.0772495406739226\\
0.113131313131313	0.0900735499423687\\
0.129292929292929	0.102874133918085\\
0.145454545454545	0.115647998522515\\
0.161616161616162	0.128391870494411\\
0.177777777777778	0.141102499905278\\
0.193939393939394	0.153776662644587\\
0.21010101010101	0.166411162871553\\
0.226262626262626	0.179002835430406\\
0.242424242424242	0.191548548226099\\
0.258585858585859	0.204045204557483\\
0.274747474747475	0.216489745405068\\
0.290909090909091	0.228879151670556\\
0.307070707070707	0.241210446365418\\
0.323232323232323	0.253480696745887\\
0.339393939393939	0.265687016391836\\
0.355555555555556	0.277826567227089\\
0.371717171717172	0.289896561478844\\
0.387878787878788	0.30189426357398\\
0.404040404040404	0.313816991970125\\
0.42020202020202	0.325662120919495\\
0.436363636363636	0.337427082163601\\
0.452525252525253	0.349109366557089\\
0.468686868686869	0.360706525619037\\
0.484848484848485	0.372216173010211\\
0.501010101010101	0.383635985934891\\
0.517171717171717	0.394963706465981\\
0.533333333333333	0.406197142792298\\
0.54949494949495	0.417334170387006\\
0.565656565656566	0.428372733096346\\
0.581818181818182	0.439310844147906\\
0.597979797979798	0.450146587077822\\
0.614141414141414	0.460878116576438\\
0.63030303030303	0.47150365925208\\
0.646464646464647	0.482021514312714\\
0.662626262626263	0.492430054165421\\
0.678787878787879	0.502727724933718\\
0.694949494949495	0.512913046892906\\
0.711111111111111	0.522984614823739\\
0.727272727272727	0.532941098284826\\
0.743434343434344	0.542781241804333\\
0.75959595959596	0.552503864991613\\
0.775757575757576	0.562107862569582\\
0.791919191919192	0.571592204328704\\
0.808080808080808	0.580955935003612\\
0.824242424242424	0.590198174073463\\
0.84040404040404	0.599318115487263\\
0.856565656565657	0.608315027315463\\
0.872727272727273	0.617188251329263\\
0.888888888888889	0.625937202509117\\
0.905050505050505	0.634561368484056\\
0.921212121212121	0.643060308903506\\
0.937373737373738	0.651433654743375\\
0.953535353535354	0.65968110754825\\
0.96969696969697	0.667802438611626\\
0.985858585858586	0.675797488096148\\
1.0020202020202	0.683666164095923\\
1.01818181818182	0.691408441643\\
1.03434343434343	0.699024361660193\\
1.05050505050505	0.706514029862441\\
1.06666666666667	0.713877615608982\\
1.08282828282828	0.721115350708616\\
1.0989898989899	0.728227528180404\\
1.11515151515152	0.73521450097215\\
1.13131313131313	0.742076680639059\\
1.14747474747475	0.748814535984982\\
1.16363636363636	0.75542859166865\\
1.17979797979798	0.761919426777358\\
1.1959595959596	0.768287673370511\\
1.21212121212121	0.774534014995493\\
1.22828282828283	0.780659185178288\\
1.24444444444444	0.786663965891287\\
1.26060606060606	0.792549186000707\\
1.27676767676768	0.79831571969601\\
1.29292929292929	0.803964484903738\\
1.30909090909091	0.809496441688096\\
1.32525252525253	0.814912590640635\\
1.34141414141414	0.820213971261342\\
1.35757575757576	0.8254016603334\\
1.37373737373737	0.830476770293846\\
1.38989898989899	0.835440447602338\\
1.40606060606061	0.840293871110162\\
1.42222222222222	0.845038250431579\\
1.43838383838384	0.849674824319579\\
1.45454545454545	0.854204859048024\\
1.47070707070707	0.858629646802131\\
1.48686868686869	0.862950504079186\\
1.5030303030303	0.867168770101306\\
1.51919191919192	0.871285805242022\\
1.53535353535354	0.875302989468395\\
1.55151515151515	0.879221720800282\\
1.56767676767677	0.883043413788345\\
1.58383838383838	0.886769498012308\\
1.6	0.890401416600884\\
};
\addplot [color=mycolor7,solid,forget plot]
  table[row sep=crcr]{%
1	0\\
1	1\\
};
\end{axis}
\end{tikzpicture}%
%
%
\definecolor{mycolor1}{rgb}{0.50000,0.50000,0.75450}%
\definecolor{mycolor2}{rgb}{0.25000,0.25000,0.63175}%
\definecolor{mycolor3}{rgb}{0.14645,0.14645,0.58091}%
\definecolor{mycolor4}{rgb}{0.06699,0.06699,0.54189}%
\definecolor{mycolor5}{rgb}{0.00000,0.00000,0.50900}%
\definecolor{mycolor6}{rgb}{0.00000,0.47170,0.46040}%
\definecolor{mycolor7}{rgb}{0.49060,0.00000,0.00000}%
\begin{tikzpicture}

\begin{axis}[%
width=0.951\figurewidth,
height=\figureheight,
at={(0\figurewidth,0\figureheight)},
scale only axis,
xmin=0,
xmax=1.6,
xtick={0,0.25,0.5,0.75,1,1.25,1.5,1.55},
xticklabels={{0},{},{0.5},{},{1},{},{},{$>1.5$}},
ymin=0,
ymax=1,
ytick={0,0.1,0.2,0.3,0.4,0.5,0.6,0.7,0.8,0.9,1},
yticklabels={\empty},
axis background/.style={fill=white},
mystyle
]
\addplot [color=mycolor1,solid,forget plot]
  table[row sep=crcr]{%
0	0.237785016286645\\
0.05	0.47557003257329\\
0.15	0.671009771986971\\
0.25	0.819543973941368\\
0.35	0.887296416938111\\
0.45	0.925732899022801\\
0.55	0.943322475570033\\
0.65	0.957654723127036\\
0.75	0.962214983713355\\
0.85	0.968078175895766\\
0.95	0.971335504885994\\
1.05	0.972638436482085\\
1.15	0.97328990228013\\
1.25	0.974592833876222\\
1.35	0.976547231270358\\
1.45	0.978501628664495\\
1.55	1\\
};
\addplot [color=mycolor2,solid,forget plot]
  table[row sep=crcr]{%
0	0.336813884577395\\
0.05	0.673627769154791\\
0.15	0.751093551573303\\
0.25	0.958374488500071\\
0.35	0.977987865105122\\
0.45	0.984196415972908\\
0.55	0.988006208550868\\
0.65	0.989558346267814\\
0.75	0.990969380555948\\
0.85	0.992239311415267\\
0.95	0.993227035416961\\
1.05	0.993791449132214\\
1.15	0.994638069705094\\
1.25	0.995061379991534\\
1.35	0.995484690277974\\
1.45	0.9957668971356\\
1.55	1\\
};
\addplot [color=mycolor3,solid,forget plot]
  table[row sep=crcr]{%
0	0.117988394584139\\
0.05	0.235976789168279\\
0.15	0.317214700193424\\
0.25	0.945454545454545\\
0.35	0.964023210831721\\
0.45	0.970212765957447\\
0.55	0.974854932301741\\
0.65	0.978336557059961\\
0.75	0.984526112185687\\
0.85	0.98568665377176\\
0.95	0.986460348162476\\
1.05	0.986847195357834\\
1.15	0.986847195357834\\
1.25	0.987620889748549\\
1.35	0.988781431334623\\
1.45	0.988781431334623\\
1.55	1\\
};
\addplot [color=mycolor4,solid,forget plot]
  table[row sep=crcr]{%
0	0.105811193629891\\
0.05	0.211622387259781\\
0.15	0.307403721001455\\
0.25	0.769083531123191\\
0.35	0.785468187734477\\
0.45	0.793048005512595\\
0.55	0.798484036444376\\
0.65	0.807595130541306\\
0.75	0.811040502258633\\
0.85	0.812265523313682\\
0.95	0.813949927264375\\
1.05	0.814792129239721\\
1.15	0.815557767399127\\
1.25	0.816323405558533\\
1.35	0.817089043717939\\
1.45	0.817471862797642\\
1.55	1\\
};
\addplot [color=mycolor5,solid,forget plot]
  table[row sep=crcr]{%
0	0.39568345323741\\
0.05	0.79136690647482\\
0.15	0.870863309352518\\
0.25	0.896043165467626\\
0.35	0.907913669064748\\
0.45	0.917985611510791\\
0.55	0.92410071942446\\
0.65	0.92841726618705\\
0.75	0.93273381294964\\
0.85	0.934892086330935\\
0.95	0.938848920863309\\
1.05	0.939928057553957\\
1.15	0.94136690647482\\
1.25	0.946043165467626\\
1.35	0.950719424460432\\
1.45	0.952877697841727\\
1.55	1\\
};
\addplot [color=mycolor6,dashed,forget plot]
  table[row sep=crcr]{%
0	0\\
0.0161616161616162	0.0128945426727563\\
0.0323232323232323	0.0257857178319707\\
0.0484848484848485	0.0386701606022311\\
0.0646464646464646	0.0515445113809416\\
0.0808080808080808	0.0644054184661257\\
0.096969696969697	0.0772495406739226\\
0.113131313131313	0.0900735499423687\\
0.129292929292929	0.102874133918085\\
0.145454545454545	0.115647998522515\\
0.161616161616162	0.128391870494411\\
0.177777777777778	0.141102499905278\\
0.193939393939394	0.153776662644587\\
0.21010101010101	0.166411162871553\\
0.226262626262626	0.179002835430406\\
0.242424242424242	0.191548548226099\\
0.258585858585859	0.204045204557483\\
0.274747474747475	0.216489745405068\\
0.290909090909091	0.228879151670556\\
0.307070707070707	0.241210446365418\\
0.323232323232323	0.253480696745887\\
0.339393939393939	0.265687016391836\\
0.355555555555556	0.277826567227089\\
0.371717171717172	0.289896561478844\\
0.387878787878788	0.30189426357398\\
0.404040404040404	0.313816991970125\\
0.42020202020202	0.325662120919495\\
0.436363636363636	0.337427082163601\\
0.452525252525253	0.349109366557089\\
0.468686868686869	0.360706525619037\\
0.484848484848485	0.372216173010211\\
0.501010101010101	0.383635985934891\\
0.517171717171717	0.394963706465981\\
0.533333333333333	0.406197142792298\\
0.54949494949495	0.417334170387006\\
0.565656565656566	0.428372733096346\\
0.581818181818182	0.439310844147906\\
0.597979797979798	0.450146587077822\\
0.614141414141414	0.460878116576438\\
0.63030303030303	0.47150365925208\\
0.646464646464647	0.482021514312714\\
0.662626262626263	0.492430054165421\\
0.678787878787879	0.502727724933718\\
0.694949494949495	0.512913046892906\\
0.711111111111111	0.522984614823739\\
0.727272727272727	0.532941098284826\\
0.743434343434344	0.542781241804333\\
0.75959595959596	0.552503864991613\\
0.775757575757576	0.562107862569582\\
0.791919191919192	0.571592204328704\\
0.808080808080808	0.580955935003612\\
0.824242424242424	0.590198174073463\\
0.84040404040404	0.599318115487263\\
0.856565656565657	0.608315027315463\\
0.872727272727273	0.617188251329263\\
0.888888888888889	0.625937202509117\\
0.905050505050505	0.634561368484056\\
0.921212121212121	0.643060308903506\\
0.937373737373738	0.651433654743375\\
0.953535353535354	0.65968110754825\\
0.96969696969697	0.667802438611626\\
0.985858585858586	0.675797488096148\\
1.0020202020202	0.683666164095923\\
1.01818181818182	0.691408441643\\
1.03434343434343	0.699024361660193\\
1.05050505050505	0.706514029862441\\
1.06666666666667	0.713877615608982\\
1.08282828282828	0.721115350708616\\
1.0989898989899	0.728227528180404\\
1.11515151515152	0.73521450097215\\
1.13131313131313	0.742076680639059\\
1.14747474747475	0.748814535984982\\
1.16363636363636	0.75542859166865\\
1.17979797979798	0.761919426777358\\
1.1959595959596	0.768287673370511\\
1.21212121212121	0.774534014995493\\
1.22828282828283	0.780659185178288\\
1.24444444444444	0.786663965891287\\
1.26060606060606	0.792549186000707\\
1.27676767676768	0.79831571969601\\
1.29292929292929	0.803964484903738\\
1.30909090909091	0.809496441688096\\
1.32525252525253	0.814912590640635\\
1.34141414141414	0.820213971261342\\
1.35757575757576	0.8254016603334\\
1.37373737373737	0.830476770293846\\
1.38989898989899	0.835440447602338\\
1.40606060606061	0.840293871110162\\
1.42222222222222	0.845038250431579\\
1.43838383838384	0.849674824319579\\
1.45454545454545	0.854204859048024\\
1.47070707070707	0.858629646802131\\
1.48686868686869	0.862950504079186\\
1.5030303030303	0.867168770101306\\
1.51919191919192	0.871285805242022\\
1.53535353535354	0.875302989468395\\
1.55151515151515	0.879221720800282\\
1.56767676767677	0.883043413788345\\
1.58383838383838	0.886769498012308\\
1.6	0.890401416600884\\
};
\addplot [color=mycolor7,solid,forget plot]
  table[row sep=crcr]{%
1	0\\
1	1\\
};
\end{axis}
\end{tikzpicture}%
%
%
\definecolor{mycolor1}{rgb}{0.50000,0.50000,0.75450}%
\definecolor{mycolor2}{rgb}{0.25000,0.25000,0.63175}%
\definecolor{mycolor3}{rgb}{0.14645,0.14645,0.58091}%
\definecolor{mycolor4}{rgb}{0.06699,0.06699,0.54189}%
\definecolor{mycolor5}{rgb}{0.00000,0.00000,0.50900}%
\definecolor{mycolor6}{rgb}{0.00000,0.47170,0.46040}%
\definecolor{mycolor7}{rgb}{0.49060,0.00000,0.00000}%
\begin{tikzpicture}

\begin{axis}[%
width=0.951\figurewidth,
height=\figureheight,
at={(0\figurewidth,0\figureheight)},
scale only axis,
xmin=0,
xmax=1.6,
xtick={0,0.25,0.5,0.75,1,1.25,1.5,1.55},
xticklabels={{0},{},{0.5},{},{1},{},{},{$>1.5$}},
ymin=0,
ymax=1,
ytick={0,0.1,0.2,0.3,0.4,0.5,0.6,0.7,0.8,0.9,1},
yticklabels={\empty},
axis background/.style={fill=white},
mystyle
]
\addplot [color=mycolor1,solid,forget plot]
  table[row sep=crcr]{%
0	0.413127413127413\\
0.05	0.826254826254826\\
0.15	0.91014391014391\\
0.25	0.935064935064935\\
0.35	0.948051948051948\\
0.45	0.956475956475957\\
0.55	0.966303966303966\\
0.65	0.96946296946297\\
0.75	0.971919971919972\\
0.85	0.973674973674974\\
0.95	0.975780975780976\\
1.05	0.978237978237978\\
1.15	0.97999297999298\\
1.25	0.982800982800983\\
1.35	0.983853983853984\\
1.45	0.984906984906985\\
1.55	1\\
};
\addplot [color=mycolor2,solid,forget plot]
  table[row sep=crcr]{%
0	0.426372498717291\\
0.05	0.852744997434582\\
0.15	0.93227296049256\\
0.25	0.963057978450487\\
0.35	0.974089276552078\\
0.45	0.978963571062083\\
0.55	0.983581323755772\\
0.65	0.987429451000513\\
0.75	0.988712160082093\\
0.85	0.990507952796306\\
0.95	0.991277578245254\\
1.05	0.993073370959466\\
1.15	0.994612621857363\\
1.25	0.995125705489995\\
1.35	0.996664956387891\\
1.45	0.996664956387891\\
1.55	1\\
};
\addplot [color=mycolor3,solid,forget plot]
  table[row sep=crcr]{%
0	0.438777262579426\\
0.05	0.877554525158853\\
0.15	0.942469517430878\\
0.25	0.964279580972008\\
0.35	0.973896616864159\\
0.45	0.980765928215696\\
0.55	0.983513652756311\\
0.65	0.986089644513138\\
0.75	0.989009101837541\\
0.85	0.990726429675425\\
0.95	0.992100291945732\\
1.05	0.99347415421604\\
1.15	0.994161085351194\\
1.25	0.995019749270136\\
1.35	0.996393611540443\\
1.45	0.99673707710802\\
1.55	1\\
};
\addplot [color=mycolor4,solid,forget plot]
  table[row sep=crcr]{%
0	0.459615941329731\\
0.05	0.919231882659462\\
0.15	0.966708482099778\\
0.25	0.981279552253209\\
0.35	0.986490398533243\\
0.45	0.989095821673261\\
0.55	0.991894239119946\\
0.65	0.993052204959954\\
0.75	0.994017176493293\\
0.85	0.995368136639969\\
0.95	0.996043616713307\\
1.05	0.996719096786645\\
1.15	0.997105085399981\\
1.25	0.997587571166651\\
1.35	0.997877062626653\\
1.45	0.998166554086654\\
1.55	1\\
};
\addplot [color=mycolor5,solid,forget plot]
  table[row sep=crcr]{%
0	0.486043042829746\\
0.05	0.972086085659493\\
0.15	0.990746096009254\\
0.25	0.99491644089982\\
0.35	0.996164500319625\\
0.45	0.997199476423853\\
0.55	0.997808285896929\\
0.65	0.998417095370004\\
0.75	0.998599738211927\\
0.85	0.998751940580195\\
0.95	0.998995464369426\\
1.05	0.999117226264041\\
1.15	0.999208547685002\\
1.25	0.999299869105963\\
1.35	0.999482511947886\\
1.45	0.999665154789809\\
1.55	1\\
};
\addplot [color=mycolor6,dashed,forget plot]
  table[row sep=crcr]{%
0	0\\
0.0161616161616162	0.0128945426727563\\
0.0323232323232323	0.0257857178319707\\
0.0484848484848485	0.0386701606022311\\
0.0646464646464646	0.0515445113809416\\
0.0808080808080808	0.0644054184661257\\
0.096969696969697	0.0772495406739226\\
0.113131313131313	0.0900735499423687\\
0.129292929292929	0.102874133918085\\
0.145454545454545	0.115647998522515\\
0.161616161616162	0.128391870494411\\
0.177777777777778	0.141102499905278\\
0.193939393939394	0.153776662644587\\
0.21010101010101	0.166411162871553\\
0.226262626262626	0.179002835430406\\
0.242424242424242	0.191548548226099\\
0.258585858585859	0.204045204557483\\
0.274747474747475	0.216489745405068\\
0.290909090909091	0.228879151670556\\
0.307070707070707	0.241210446365418\\
0.323232323232323	0.253480696745887\\
0.339393939393939	0.265687016391836\\
0.355555555555556	0.277826567227089\\
0.371717171717172	0.289896561478844\\
0.387878787878788	0.30189426357398\\
0.404040404040404	0.313816991970125\\
0.42020202020202	0.325662120919495\\
0.436363636363636	0.337427082163601\\
0.452525252525253	0.349109366557089\\
0.468686868686869	0.360706525619037\\
0.484848484848485	0.372216173010211\\
0.501010101010101	0.383635985934891\\
0.517171717171717	0.394963706465981\\
0.533333333333333	0.406197142792298\\
0.54949494949495	0.417334170387006\\
0.565656565656566	0.428372733096346\\
0.581818181818182	0.439310844147906\\
0.597979797979798	0.450146587077822\\
0.614141414141414	0.460878116576438\\
0.63030303030303	0.47150365925208\\
0.646464646464647	0.482021514312714\\
0.662626262626263	0.492430054165421\\
0.678787878787879	0.502727724933718\\
0.694949494949495	0.512913046892906\\
0.711111111111111	0.522984614823739\\
0.727272727272727	0.532941098284826\\
0.743434343434344	0.542781241804333\\
0.75959595959596	0.552503864991613\\
0.775757575757576	0.562107862569582\\
0.791919191919192	0.571592204328704\\
0.808080808080808	0.580955935003612\\
0.824242424242424	0.590198174073463\\
0.84040404040404	0.599318115487263\\
0.856565656565657	0.608315027315463\\
0.872727272727273	0.617188251329263\\
0.888888888888889	0.625937202509117\\
0.905050505050505	0.634561368484056\\
0.921212121212121	0.643060308903506\\
0.937373737373738	0.651433654743375\\
0.953535353535354	0.65968110754825\\
0.96969696969697	0.667802438611626\\
0.985858585858586	0.675797488096148\\
1.0020202020202	0.683666164095923\\
1.01818181818182	0.691408441643\\
1.03434343434343	0.699024361660193\\
1.05050505050505	0.706514029862441\\
1.06666666666667	0.713877615608982\\
1.08282828282828	0.721115350708616\\
1.0989898989899	0.728227528180404\\
1.11515151515152	0.73521450097215\\
1.13131313131313	0.742076680639059\\
1.14747474747475	0.748814535984982\\
1.16363636363636	0.75542859166865\\
1.17979797979798	0.761919426777358\\
1.1959595959596	0.768287673370511\\
1.21212121212121	0.774534014995493\\
1.22828282828283	0.780659185178288\\
1.24444444444444	0.786663965891287\\
1.26060606060606	0.792549186000707\\
1.27676767676768	0.79831571969601\\
1.29292929292929	0.803964484903738\\
1.30909090909091	0.809496441688096\\
1.32525252525253	0.814912590640635\\
1.34141414141414	0.820213971261342\\
1.35757575757576	0.8254016603334\\
1.37373737373737	0.830476770293846\\
1.38989898989899	0.835440447602338\\
1.40606060606061	0.840293871110162\\
1.42222222222222	0.845038250431579\\
1.43838383838384	0.849674824319579\\
1.45454545454545	0.854204859048024\\
1.47070707070707	0.858629646802131\\
1.48686868686869	0.862950504079186\\
1.5030303030303	0.867168770101306\\
1.51919191919192	0.871285805242022\\
1.53535353535354	0.875302989468395\\
1.55151515151515	0.879221720800282\\
1.56767676767677	0.883043413788345\\
1.58383838383838	0.886769498012308\\
1.6	0.890401416600884\\
};
\addplot [color=mycolor7,solid,forget plot]
  table[row sep=crcr]{%
1	0\\
1	1\\
};
\end{axis}
\end{tikzpicture}%
%
%
\definecolor{mycolor1}{rgb}{0.50000,0.50000,0.75450}%
\definecolor{mycolor2}{rgb}{0.25000,0.25000,0.63175}%
\definecolor{mycolor3}{rgb}{0.14645,0.14645,0.58091}%
\definecolor{mycolor4}{rgb}{0.06699,0.06699,0.54189}%
\definecolor{mycolor5}{rgb}{0.00000,0.00000,0.50900}%
\definecolor{mycolor6}{rgb}{0.00000,0.47170,0.46040}%
\definecolor{mycolor7}{rgb}{0.49060,0.00000,0.00000}%
\begin{tikzpicture}

\begin{axis}[%
width=0.951\figurewidth,
height=\figureheight,
at={(0\figurewidth,0\figureheight)},
scale only axis,
xmin=0,
xmax=1.6,
xtick={0,0.25,0.5,0.75,1,1.25,1.5,1.55},
xticklabels={{0},{},{0.5},{},{1},{},{},{$>1.5$}},
ymin=0,
ymax=1,
ytick={0,0.1,0.2,0.3,0.4,0.5,0.6,0.7,0.8,0.9,1},
yticklabels={\empty},
axis background/.style={fill=white},
mystyle
]
\addplot [color=mycolor1,solid,forget plot]
  table[row sep=crcr]{%
0	0.315881326352531\\
0.05	0.631762652705061\\
0.15	0.778359511343805\\
0.25	0.842931937172775\\
0.35	0.884816753926702\\
0.45	0.910994764397906\\
0.55	0.921465968586387\\
0.65	0.93717277486911\\
0.75	0.945898778359511\\
0.85	0.952879581151832\\
0.95	0.963350785340314\\
1.05	0.965095986038394\\
1.15	0.966841186736475\\
1.25	0.970331588132635\\
1.35	0.975567190226876\\
1.45	0.984293193717278\\
1.55	1\\
};
\addplot [color=mycolor2,solid,forget plot]
  table[row sep=crcr]{%
0	0.45137464095199\\
0.05	0.90274928190398\\
0.15	0.959171112022979\\
0.25	0.974148543290931\\
0.35	0.981329503487895\\
0.45	0.985843249897415\\
0.55	0.988920804267542\\
0.65	0.990767336889618\\
0.75	0.991793188346327\\
0.85	0.992819039803037\\
0.95	0.994050061551087\\
1.05	0.994870742716455\\
1.15	0.995691423881822\\
1.25	0.997127615921215\\
1.35	0.99774312679524\\
1.45	0.998974148543291\\
1.55	1\\
};
\addplot [color=mycolor3,solid,forget plot]
  table[row sep=crcr]{%
0	0.456640134059489\\
0.05	0.913280268118978\\
0.15	0.961876832844575\\
0.25	0.973397570171764\\
0.35	0.981357352325094\\
0.45	0.98533724340176\\
0.55	0.988060326770004\\
0.65	0.990783410138249\\
0.75	0.99266862170088\\
0.85	0.993087557603687\\
0.95	0.994134897360704\\
1.05	0.994763301214914\\
1.15	0.996229576874738\\
1.25	0.996648512777545\\
1.35	0.997905320485966\\
1.45	0.998324256388773\\
1.55	1\\
};
\addplot [color=mycolor4,solid,forget plot]
  table[row sep=crcr]{%
0	0.23943661971831\\
0.05	0.47887323943662\\
0.15	0.690140845070423\\
0.25	0.774647887323944\\
0.35	0.830985915492958\\
0.45	0.845070422535211\\
0.55	0.873239436619718\\
0.65	0.873239436619718\\
0.75	0.873239436619718\\
0.85	0.873239436619718\\
0.95	0.873239436619718\\
1.05	0.887323943661972\\
1.15	0.887323943661972\\
1.25	0.901408450704225\\
1.35	0.901408450704225\\
1.45	0.901408450704225\\
1.55	1\\
};
\addplot [color=mycolor5,solid,forget plot]
  table[row sep=crcr]{%
0	0.415841584158416\\
0.05	0.831683168316832\\
0.15	0.935643564356436\\
0.25	0.96039603960396\\
0.35	0.965346534653465\\
0.45	0.965346534653465\\
0.55	0.965346534653465\\
0.65	0.965346534653465\\
0.75	0.965346534653465\\
0.85	0.965346534653465\\
0.95	0.965346534653465\\
1.05	0.965346534653465\\
1.15	0.965346534653465\\
1.25	0.965346534653465\\
1.35	0.965346534653465\\
1.45	0.965346534653465\\
1.55	1\\
};
\addplot [color=mycolor6,dashed,forget plot]
  table[row sep=crcr]{%
0	0\\
0.0161616161616162	0.0128945426727563\\
0.0323232323232323	0.0257857178319707\\
0.0484848484848485	0.0386701606022311\\
0.0646464646464646	0.0515445113809416\\
0.0808080808080808	0.0644054184661257\\
0.096969696969697	0.0772495406739226\\
0.113131313131313	0.0900735499423687\\
0.129292929292929	0.102874133918085\\
0.145454545454545	0.115647998522515\\
0.161616161616162	0.128391870494411\\
0.177777777777778	0.141102499905278\\
0.193939393939394	0.153776662644587\\
0.21010101010101	0.166411162871553\\
0.226262626262626	0.179002835430406\\
0.242424242424242	0.191548548226099\\
0.258585858585859	0.204045204557483\\
0.274747474747475	0.216489745405068\\
0.290909090909091	0.228879151670556\\
0.307070707070707	0.241210446365418\\
0.323232323232323	0.253480696745887\\
0.339393939393939	0.265687016391836\\
0.355555555555556	0.277826567227089\\
0.371717171717172	0.289896561478844\\
0.387878787878788	0.30189426357398\\
0.404040404040404	0.313816991970125\\
0.42020202020202	0.325662120919495\\
0.436363636363636	0.337427082163601\\
0.452525252525253	0.349109366557089\\
0.468686868686869	0.360706525619037\\
0.484848484848485	0.372216173010211\\
0.501010101010101	0.383635985934891\\
0.517171717171717	0.394963706465981\\
0.533333333333333	0.406197142792298\\
0.54949494949495	0.417334170387006\\
0.565656565656566	0.428372733096346\\
0.581818181818182	0.439310844147906\\
0.597979797979798	0.450146587077822\\
0.614141414141414	0.460878116576438\\
0.63030303030303	0.47150365925208\\
0.646464646464647	0.482021514312714\\
0.662626262626263	0.492430054165421\\
0.678787878787879	0.502727724933718\\
0.694949494949495	0.512913046892906\\
0.711111111111111	0.522984614823739\\
0.727272727272727	0.532941098284826\\
0.743434343434344	0.542781241804333\\
0.75959595959596	0.552503864991613\\
0.775757575757576	0.562107862569582\\
0.791919191919192	0.571592204328704\\
0.808080808080808	0.580955935003612\\
0.824242424242424	0.590198174073463\\
0.84040404040404	0.599318115487263\\
0.856565656565657	0.608315027315463\\
0.872727272727273	0.617188251329263\\
0.888888888888889	0.625937202509117\\
0.905050505050505	0.634561368484056\\
0.921212121212121	0.643060308903506\\
0.937373737373738	0.651433654743375\\
0.953535353535354	0.65968110754825\\
0.96969696969697	0.667802438611626\\
0.985858585858586	0.675797488096148\\
1.0020202020202	0.683666164095923\\
1.01818181818182	0.691408441643\\
1.03434343434343	0.699024361660193\\
1.05050505050505	0.706514029862441\\
1.06666666666667	0.713877615608982\\
1.08282828282828	0.721115350708616\\
1.0989898989899	0.728227528180404\\
1.11515151515152	0.73521450097215\\
1.13131313131313	0.742076680639059\\
1.14747474747475	0.748814535984982\\
1.16363636363636	0.75542859166865\\
1.17979797979798	0.761919426777358\\
1.1959595959596	0.768287673370511\\
1.21212121212121	0.774534014995493\\
1.22828282828283	0.780659185178288\\
1.24444444444444	0.786663965891287\\
1.26060606060606	0.792549186000707\\
1.27676767676768	0.79831571969601\\
1.29292929292929	0.803964484903738\\
1.30909090909091	0.809496441688096\\
1.32525252525253	0.814912590640635\\
1.34141414141414	0.820213971261342\\
1.35757575757576	0.8254016603334\\
1.37373737373737	0.830476770293846\\
1.38989898989899	0.835440447602338\\
1.40606060606061	0.840293871110162\\
1.42222222222222	0.845038250431579\\
1.43838383838384	0.849674824319579\\
1.45454545454545	0.854204859048024\\
1.47070707070707	0.858629646802131\\
1.48686868686869	0.862950504079186\\
1.5030303030303	0.867168770101306\\
1.51919191919192	0.871285805242022\\
1.53535353535354	0.875302989468395\\
1.55151515151515	0.879221720800282\\
1.56767676767677	0.883043413788345\\
1.58383838383838	0.886769498012308\\
1.6	0.890401416600884\\
};
\addplot [color=mycolor7,solid,forget plot]
  table[row sep=crcr]{%
1	0\\
1	1\\
};
\end{axis}
\end{tikzpicture}%
  \caption{Empirical cummulative distriubtion function (CDF)
  		   of true local errors $\xi_n$ divided by
           the estimated local errors $(Q(t_n))_{00}^{1/2}$.
           Ticks on the $y$-axis are spaced at $0.1$ intervals
           from $0$ to $1$.
           Values less than $1$ (red line) are \emph{over-estimated}
           leading to a conservative step size adaptation.
           Green dashed line shows the CDF of the $\chi(1)$-distribution
           which implies that the empirical distribution has weaker tails.
           See text for more details.}
  \label{fig:errorEstimation}
\end{figure*}
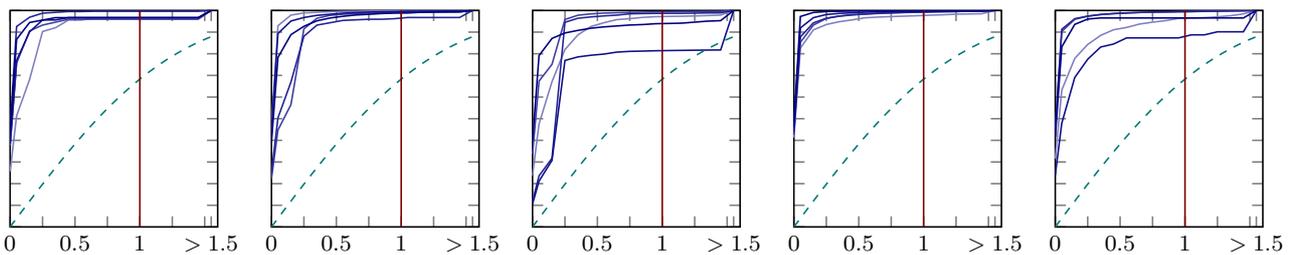

\edit{
In addition to the benchmark results,
we analyze the error estimation model from
a probabilistic perspective.
Figure~\ref{fig:errorEstimation} plots the
cumulative distribution function (CDF) of
the local error $\xi_n$, as defined above,
divided by the estimated local error
$(Q(\tn))_{00}^{1/2} = (\ssq_n \bar{Q}(h_n))_{00}^{1/2}$
for each set of five tasks
(different blue colored lines) of each
of the five problem classes (figures from left to right).
Under the algorithm's internal model, the error is assumed to be Gaussian distributed:
\begin{equation}
	P(y_\tn\g\hat{y}_\tn) = \N(y_\tn; \hat{y}_\tn, (\mat{Q}(h_n))_{00}),
\end{equation}
Hence, if that model were a perfect fit, the scaled absolute error plotted in this figure would be $\chi$-distributed:
\begin{equation}
	P(\abs{y_\tn - \hat{y}\tn}(\mat{Q}(h_n))_{00}^{-1/2}) = \chi(1).
\end{equation}
The comparison with the CDF of $\chi(1)$ shows
that the true local error has weaker tails than
the predicted $\chi$-distribution.
}

\edit{
So, as intended, the error estimator is typically a 
conservative one.
}


While the probabilistic method does not achieve the same
high performance as modern higher order codes, the performance
matches the results of a production Runge--Kutta code of the
same order. This is of particular interest since applications
in the low accuracy regime could benefit the most from accurate
error indicators \cite{gear1981anything}.


\section{Conclusions}
\label{sec:conclusions}

We proposed a probabilistic inference model for the numerical
solution of ODEs and showed the connections with established
methods. In particular, we showed how probabilistic inference in
Gauss--Markov systems given by a linear time-invariant stochastic differential
equations leads to Nordsieck-type methods. The
maximum a posteriori estimate of the once integrated Wiener process IWP$(1)$
is equivalent to the trapezoidal rule.
The twice integrated Wiener process IWP$(2)$ is equivalent to a third order
Nordsieck-type method which can be thought of as a spline-based multistep method.
We demonstrated the practicality of this probabilistic IVP
solver by comparing against other state-of-the-art implementations.



The probabilistic formulation has already proven to be beneficial in
larger chains of computations involving boundary value problems
\cite{schober2014probabilistic,hauberg2015random}.
While the method presented in this paper is restricted to IVPs, there
has also been work on extending the formalism of splines to boundary
value problems \cite{Mazzia2006BSpline,Mazzia2009continuous}.
We expect that similar classical guarantees should be transferable to
probabilistic \editone{boundary value problem} solvers as well. Conversely, the probabilistic
treatment of the IVP may be beneficial in bigger pipelines as well (cf.~\cite{o.13:_bayes_uncer_quant_differ_equat}).

\appendix

\section{Detailed results}
\label{sec:detailed-results}
\begin{figure*}
  \centering
  \scriptsize{
  \setlength{\figurewidth}{0.9\textwidth}
  \setlength{\figureheight}{0.19\textheight}
  \input{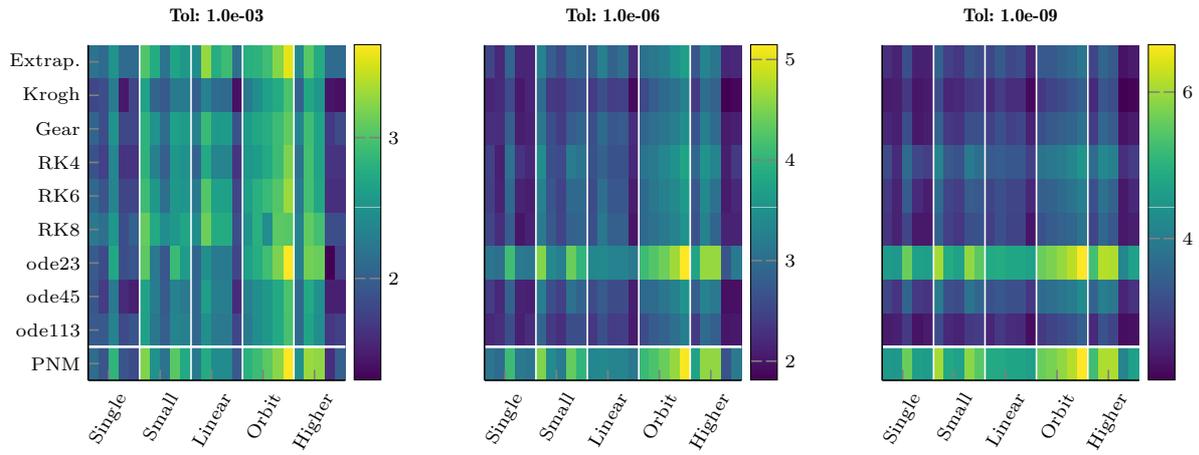}}
  \caption{$\log_{10}$(\#FE), the number of function
    evaluations in logarithmic scale, for all tested methods
    and individual problems.}
  \label{fig:fcnCalls}
\end{figure*}

\begin{figure*}
  \centering
  \scriptsize{
  \setlength{\figurewidth}{0.9\textwidth}
  \setlength{\figureheight}{0.19\textheight}
  \input{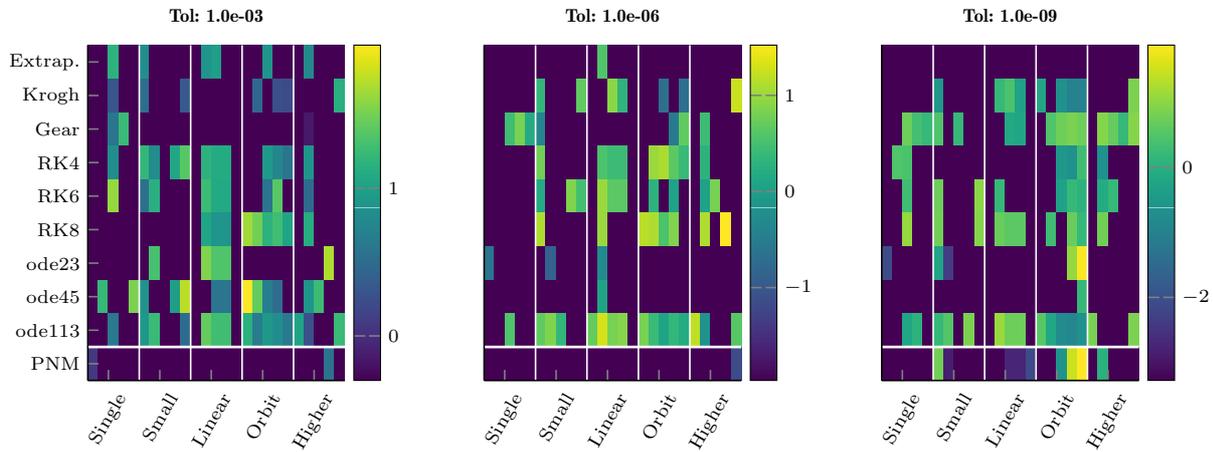}}
  \caption{$\log_{10}(|\{\xi_n\g \xi_n > \tol, n = 1, \dotsc, N\}|N^{-1})$,
    the percent of deceived steps in logarithmic scale, for all
    tested methods and individual problems.}
  \label{fig:percDecv}
\end{figure*}

\begin{figure*}
  \centering
  \scriptsize{
  \setlength{\figurewidth}{0.9\textwidth}
  \setlength{\figureheight}{0.19\textheight}
  \input{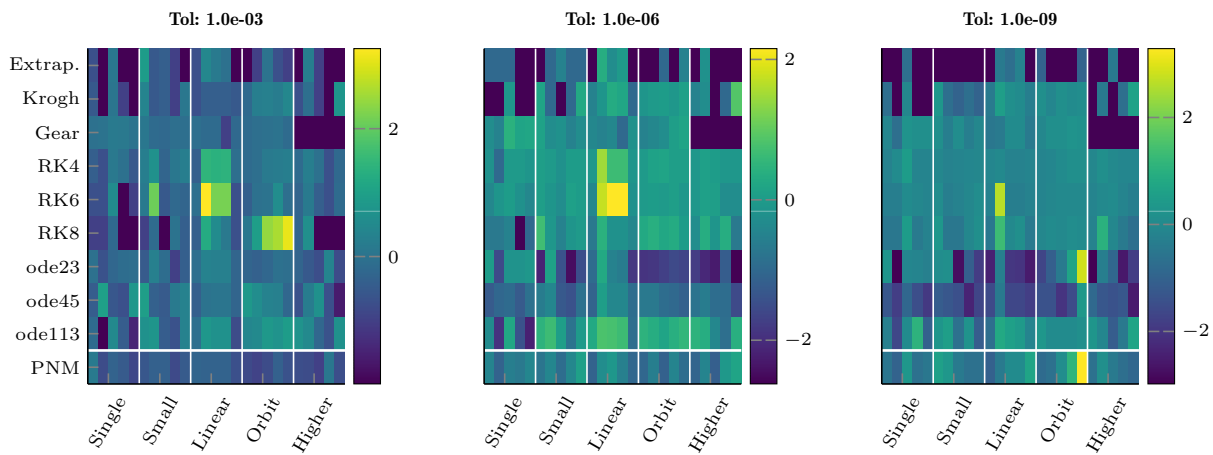}}
  \caption{$\log_{10}(\max \{\xi_n[h_n \tol]^{-1}\g n = 1, \dotsc, N\})$,
    the maximum error per unit step in logarithmic scale, for all 
    tested methods and individual problems}
  \label{fig:maxError}
\end{figure*}

Figures~\ref{fig:fcnCalls}, \ref{fig:percDecv} and \ref{fig:maxError}
in this section present detailed results from
the \texttt{DETEST} test set. For a detailed
description see Sec.~\ref{sec:experiments} and
\cite{hull1972comparing}.

\section{A fourth order four step Runge--Kutta formula expressed as
LTI SDE filtering solution}
\label{sec:rk-derivations}

Runge--Kutta (RK) formulas are a family of one step solvers
for ODEs. At time $\tnpo$, the numerical approximation is
defined
\begin{equation}
\begin{aligned}
	y_\tnpo &= y_\tn + h \sum_{i=1}^s b_i k_{i,n}\label{eq:runge-kutta}\\
	k_{i,n} &= f(\tn + h c_i,\,y_\tn + h \sum_{j=1}^{s} a_{ij} k_{j,n})
\end{aligned}
\end{equation}
The parameters $a_{ij}, c_i$ and $b_i$ are usually expressed as
a matrix $\mat{A}$ and two vectors $\vec{b}, \vec{c}$, written
compactly in a so-called Butcher-tableau:
\begin{equation}
\begin{array}{c|ccc}
c_1 & a_{11} & \hdots & a_{1s}\\
\vdots & \vdots & & \vdots\\
c_s & a_{s1} & \hdots & a_{ss}\\
\hline
& b_1 & \hdots & b_s
\end{array}
\end{equation}
If the matrix $\mat{A}$ is strictly lower triangular and $c_1 = 0$,
Equation~\eqref{eq:runge-kutta} simplifies to an iterative procedure of
explicit equations
\begin{equation}
	\begin{aligned}
		k_{1,n} &= f(\tn,\,y_\tn),\\
		k_{i,n} &= f(\tn + c_i h,\, y_\tn + h \sum_{j=1}^{i-1} a_{ij} k_{j,n}),\q i = 2, \dotsc, s\\
		y_\tnpo &= y_\tn + h \sum_{i=1}^s b_i k_{i,n},\label{eq:explicit-RK}
	\end{aligned}
\end{equation}
and in this case, the formula is called an \emph{explicit RK method}.
A RK method is said to be of order $q$, iff
\begin{equation}
	\abs{y(t_0 + h) - y_{t_1}} \leq Kh^{q+1},
\end{equation}
where $K$ is a constant independent of $h$. In this case, the
global error is of order $h^q$. It can be shown that there
exist RK methods for which the order $q = s$ the number of
vector field evaluations up to and including $q = s = 4$.
Furthermore, there are no RK methods for which $q = s$, if
$q \geq 5$. For a full description, we refer the reader
to \cite{hairer87:_solvin_ordin_differ_equat_i}.

Schober et al.~\cite{schober2014nips} have presented
probabilistic models whose MAP at $t_0 + h = t_1$ is
equivalent to an explicit RK method of type Equation~\eqref{eq:explicit-RK}
in the case of $q \in \{1, 2, 3\}$. The corresponding
probabilistic model is given by
\begin{equation}
	\diff \vec{X} = \mat{U}_{q+1} \vec{X} \diff t + \vec{e}_q \diff W,\label{eq:iwp-apx}
\end{equation}
where $(\mat{U}_{q+1})_{ij} = \delta_{i,j+1},\,i,j=0, \dotsc, q$ is the
$q+1$-dimensional upper shift matrix and $\vec{e}_q = (0, \dotsc, 0, 1)\Trans$
is the $q+1$\textsuperscript{th} standard basis vector.
The process~\eqref{eq:iwp-apx} is known as
\emph{$q$-times integrated Wiener process} IWP($q$).
The initial conditions are
\begin{equation}
\begin{aligned}
	\vec{m}_\tmo &= \vec{0}, & \mat{C}_\tmo &= \lim_{\mathcal{H} \to \infty} \mat{Q}(\mathcal{H}).
\end{aligned}
\end{equation}
See \cite{schober2014nips} for details.

Assume that Algorithm ~\ref{alg:prob-nordsieck} has been
run until the loop has been evaluated four times with
$h_n$ such that $(t_0, \dots, t_3) = \vec{c}\Trans h = (0, uh, vh, 1)$,
where $u, v$ are two constant in $[0, 1]$ chosen by the user.
Then, $lim_{\mathcal{H} \to \infty} \vec{m}_{t_3}$ for the IWP($4$) is
\begin{align*}
	(\vec{m}_{t_3})_0 &= y_0 + h \frac{1 - 2(u+v) + 6uv}{12uv} z_0
	                   + h \frac{2v-1}{12u(u-v)(u-1)} z_1\\
	            &\phantom{= y_0 } + h \frac{1 - 2u}{12v(u-v)(v-1)} z_2
	                   + h \frac{3-4(u+v)+6uv}{12(u-1)(v-1)} z_3\\
	(\vec{m}_{t_3})_1 &= z_3\\
	(\vec{m}_{t_3})_2 &= \frac{1}{h}\frac{u+v-uv-1}{uv} z_0
	               + \frac{1}{h}\frac{1-v}{u(u-v)(u-1)} z_1\\
	            &+ \frac{1}{h}\frac{u-1}{v(u-v)(v-1)} z_2
	               + \frac{1}{h}\frac{3-2(u+v)+uv}{(u-1)(v-1)} z_3\\
	(\vec{m}_{t_3})_3 &= \frac{1}{h^2}\frac{2(u+v-2)}{uv} z_0
	               + \frac{1}{h^2}\frac{2(2-v)}{u(u-v)(u-1)} z_1\\
	            &+ \frac{1}{h^2}\frac{2(u-2)}{v(u-v)(v-1)} z_2
	               + \frac{1}{h^2}\frac{2(3-u-v)}{(u-1)(v-1)} z_3\\
	(\vec{m}_{t_3})_4 &= \frac{1}{h^3}\frac{-6}{uv} z_0
	               + \frac{1}{h^3}\frac{6}{u(u-v)(u-1)} z_1\\
	            &+ \frac{1}{h^3}\frac{-6}{v(u-v)(v-1)} z_2
	               + \frac{1}{h^3}\frac{6}{(u-1)(v-1)} z_3	 	                             
\end{align*}
Furthermore, we get the following algebraic equations for
the elements of the covariance matrix
$lim_{\mathcal{H} \to \infty} \mat{C}_{t_3}$:
\begin{align*}
	(\mat{C}_{t_3})_{00} &= \ssq h^9 [6u^6v^2 - 3u^6v + 6u^5v^3 - 27u^5v^2 + 20u^5v\\
	  &\phantom{= \ssq h^9}  - 4u^5 + 6u^4v^4 - 27u^4v^3 + 28u^4v^2 - 12u^4v\\
	  &\phantom{= \ssq h^9}   + 2u^4 + 6u^3v^5 - 27u^3v^4+ 28u^3v^3 - 12u^3v^2\\
	  &\phantom{= \ssq h^9}   + 2u^3v + 6u^2v^6- 27u^2v^5 + 28u^2v^4 + 68u^2v^3\\
	  &\phantom{= \ssq h^9}   - 78u^2v^2 + 20u^2v - 9uv^6+ 38uv^5 - 42uv^4\\
	  &\phantom{= \ssq h^9}   - 48uv^3 + 70uv^2 - 20uv + 3v^6 - 13v^5 + 17v^4\\
	  &\phantom{= \ssq h^9}    + 5v^3 - 15v^2+ 5v]/[725760v(1 - u)]\\
\end{align*}
\begin{align*}
	(\mat{C}_{t_3})_{01} &= 0\\
	(\mat{C}_{t_3})_{02} &= \ssq h^7 [(v - 1)(3u^6v + 3u^5v^2 - 16u^5v + 6u^5 + 3u^4v^3\\
	          &\phantom{= \ssq h^7}  - 16u^4v^2 + 14u^4v - 4u^4 + 3u^3v^4 - 16u^3v^3\\
	          &\phantom{= \ssq h^7}  + 14u^3v^2 - 4u^3v + 3u^2v^5 - 16u^2v^4 + 14u^2v^3\\
	          &\phantom{= \ssq h^7}  + 76u^2v^2 - 40u^2v - 6uv^5 + 29uv^4 - 24uv^3\\
	          &\phantom{= \ssq h^7}  - 85uv^2 + 50uv + 3v^5 - 14v^4 + 15v^3 + 20v^2\\
	          &\phantom{= \ssq h^7} - 15v)]/[120960v(u - 1)]\\
	(\mat{C}_{t_3})_{03} &= \ssq h^6 [3u^6v^2 - 6u^6v + 3u^5v^3 - 18u^5v^2  + 32u^5v\\
	          &\phantom{= \ssq h^6} - 10u^5 + 3u^4v^4 - 18u^4v^3 + 40u^4v^2 - 30u^4v\\
	          &\phantom{= \ssq h^6} + 8u^4 + 3u^3v^5 - 18u^3v^4 + 40u^3v^3 - 30u^3v^2\\
	          &\phantom{= \ssq h^6} + 8u^3v + 3u^2v^6 - 18u^2v^5 + 40u^2v^4 + 50u^2v^3\\
	          &\phantom{= \ssq h^6} - 192u^2v^2 + 80u^2v - 9uv^6 + 41uv^5 - 69uv^4\\
	          &\phantom{= \ssq h^6} - 81uv^3 + 271uv^2 - 117uv + 6v^6 - 28v^5 \\
	          &\phantom{= \ssq h^6} + 50v^4 + 2v^3 - 78v^2 + 39v]/[60480v(u - 1)]\\
\end{align*}
\begin{align*}
	(\mat{C}_{t_3})_{04} &= \ssq h^5 [3u^6v + 3u^5v^2 - 12u^5v + 4u^5 + 3u^4v^3\\
	         &\phantom{= \ssq h^5}   - 12u^4v^2 + 12u^4v - 4u^4 + 3u^3v^4 - 12u^3v^3\\
	         &\phantom{= \ssq h^5}   + 12u^3v^2- 4u^3v + 3u^2v^5 - 12u^2v^4 + 12u^2v^3\\
	         &\phantom{= \ssq h^5}   + 76u^2v^2 - 40u^2v+ 3uv^6 - 12uv^5 + 12uv^4\\
	         &\phantom{= \ssq h^5}   + 16uv^3 - 140uv^2 + 72uv - 3v^6 + 13v^5 - 19v^4\\
	         &\phantom{= \ssq h^5}    + 5v^3 + 45v^2 - 27v]/[20160v(1 - u)]\\
	(\mat{C}_{t_3})_{11} &= 0\\
	(\mat{C}_{t_3})_{12} &= 0\\
	(\mat{C}_{t_3})_{13} &= 0\\
	(\mat{C}_{t_3})_{14} &= 0.\\
\end{align*}
The last four equations are a consequence of the noise-free observation $z_3$
at $t_0 + c_4h = t_0 + h = t_3$. The remaining entries are given by the expressions
\begin{align*}
	(\mat{C}_{t_3})_{22} &= \ssq h^5 [(v - 1)^2(u^4 + u^3v + u^2v^2 + uv^3 - 10uv\\
	  &\phantom{= \ssq h^5} - 2v^3 + 2v^2 + 6v)]/[2520v] \\
	(\mat{C}_{t_3})_{23} &= \ssq h^4 [(v - 1)(u^5v - 3u^5 + u^4v^2 - 5u^4v + 4u^4\\
	           &\phantom{= \ssq h^4} + u^3v^3 - 5u^3v^2 + 4u^3v + u^2v^4 - 5u^2v^3\\
	           &\phantom{= \ssq h^4} - 16u^2v^2 + 40u^2v - 5uv^4 + 15uv^3 + 37uv^2\\
	           &\phantom{= \ssq h^4} - 77uv + 5v^4 - 15v^3 - 11v^2 + 33v)]/[2520v(u - 1)]\\
	(\mat{C}_{t_3})_{24} &= \ssq h^3 [(v - 1)(- u^5 - u^4v + 2u^4 - u^3v^2 + 2u^3v\\
	          &\phantom{= \ssq h^3}  - u^2v^3 + 2u^2v^2 + 20u^2v - uv^4 + 2uv^3\\
	          &\phantom{= \ssq h^3} + 5uv^2 - 50uv + 2v^4 - 5v^3 + 25v)]/[840v(u - 1)]\\
	(\mat{C}_{t_3})_{33} &= \ssq h^3 [u^5v - 2u^5 + 2u^4v^2 - 6u^4v + 4u^4 + 2u^3v^3\\
	          &\phantom{= \ssq h^3}  - 6u^3v^2 + 4u^3v + 2u^2v^4 + 4u^2v^3 - 36u^2v^2\\
	          &\phantom{= \ssq h^3}   + 40u^2v + uv^5 - 12uv^4 - 12uv^3 + 104uv^2 - 96uv\\
	          &\phantom{= \ssq h^3}  - 2v^5 + 16v^4 - 8v^3 - 48v^2 + 48v]/[630v(1 - u)]\\
	(\mat{C}_{t_3})_{34} &= \ssq h^2[u^5 + 3u^4v - 4u^4 + 3u^3v^2 - 4u^3v + 3u^2v^3\\
      &\phantom{= \ssq h^2} + 16u^2v^2 - 40u^2v + 3uv^4 + uv^3 - 95uv^2 + 135uv\\
      &\phantom{= \ssq h^2} + v^5- 10v^4 + 14v^3 + 54v^2 - 81v]/[420v(u - 1)]\\
	(\mat{C}_{t_3})_{44} &= \ssq h [u^4 + u^3v + u^2v^2 + 10u^2v + uv^3 + 20uv^2 \\
	 &\phantom{= \ssq h} - 60uv + v^4 - 5v^3 - 15v^2 + 45v]/[70v(1-u)]
\end{align*}
which defines the entire matrix since $\mat{C}_{t_3}$ is a
symmetric matrix.

The proof that this specific choice of $\vec{c}$ is analogous to the
proofs given in \cite{schober2014nips} and can be checked with laborious
algebra.

We would like to point to one more detail: although it can
easily be checked that $(\vec{m}_{t_3})_0$ is of the required form
to produce the RK prediction, this does not suffice to show that
this choice of evaluation knots produces Runge--Kutta.
The space of suitable parameters to produce a $q$th-order method
is constrained by the expansion to match the Taylor coefficients.
In the case of the IWP($q$), where each subsequent evaluation increases
the order of the polynomial approximation, this entails that each
partial RK-step needs to be a RK method of its own right to produce
an overall RK method of high order. One can think about this as a
bigger set of constraints that needs to be fulfilled. As a consequence,
this also entails that there is no meaningful interpretation of
RK methods with $v \neq \nicefrac{2}{3}$ in the case of the IWP($3$)
as has erroneously been conjectured in \cite{schober2014nips}.

For complete details, see \cite{hairer87:_solvin_ordin_differ_equat_i,schober2014nips}.

\begin{acknowledgements}
\editone{The authors thank Hans Kersting for valuable discussions and helpful comments on the manuscript.
The authors also thank the feedback of the anonymous reviewers which helped to improve the
presentation significantly.}
\end{acknowledgements}

\bibliographystyle{spmpsci}      
\bibliography{bibfile}   

\end{document}